\documentclass{amsart}
\usepackage{graphicx}
\usepackage{amssymb}
\usepackage{color}
\usepackage[dvipsnames]{xcolor}
\usepackage{tikz}
\usepackage{hyperref}
\usepackage{cleveref}
\usepackage{float}
\usepackage{mathtools}
\usepackage{enumitem}
\usepackage{cancel}
\usepackage{ifthen}

\title{The Hurwitz problem for abelian differentials}
\date{October 2025}

\DeclareMathOperator{\ord}{ord}
\DeclareMathOperator{\Homeo}{Homeo}
\DeclareMathOperator{\Aff}{Aff}

\DeclareMathOperator{\trans}{Trans}
\DeclareMathOperator{\Aut}{Aut}
\DeclareMathOperator{\Inn}{Inn}
\DeclareMathOperator{\GL}{GL}
\DeclareMathOperator{\SL}{SL}
\DeclareMathOperator{\PSL}{PSL}
\DeclareMathOperator{\SO}{SO}
\DeclareMathOperator{\Id}{Id}
\DeclareMathOperator{\tr}{tr}
\DeclareMathOperator{\diag}{diag}
\DeclareMathOperator{\Gal}{Gal}

\newcommand{\qi}{\mathbf{i}}
\newcommand{\qj}{\mathbf{j}}
\newcommand{\qk}{\mathbf{k}}

\newcommand{\divides}{\mid}
\newcommand{\ndivides}{\nmid}
\newcommand{\translations}{\mathsf{t}}
\newcommand{\orderm}{\mathsf{m}}
\newcommand{\constant}{\mathsf{c}}
\newcommand{\genera}{\mathsf{G}}
\newcommand{\st}{\ \mid\ }

\newcommand{\Dic}{\mathrm{Dic}}
\newcommand{\D}{\mathrm{D}}
\newcommand{\SD}{\mathrm{SD}}
\newcommand{\M}{\mathrm{M}}
\newcommand{\Sym}{\mathrm{S}}
\newcommand{\A}{\mathrm{A}}
\newcommand{\Q}{\mathrm{Q}}

\newcommand{\N}{\mathbf{N}}
\newcommand{\C}{\mathbf{C}}
\newcommand{\Z}{\mathbf{Z}}

\newcommand{\FF}{\mathbb{F}}
\newcommand{\NN}{\mathbb{N}}
\newcommand{\RR}{\mathbb{R}}
\newcommand{\TT}{\mathbb{T}}
\newcommand{\ZZ}{\mathbb{Z}}

\newcommand{\calH}{\mathcal{H}}

\numberwithin{equation}{section}

\numberwithin{figure}{section}

\newtheorem{Thm}[equation]{Theorem}
\newtheorem{Lem}[equation]{Lemma}
\newtheorem{Prop}[equation]{Proposition}
\newtheorem{Rema}[equation]{Remark}

\newtheorem{Cor}[equation]{Corollary}

\newtheorem*{Notation}{Notation}

\newtheorem{ThmA}{Theorem}

\crefname{ThmA}{Theorem}{Theorems}
\crefname{Thm}{Theorem}{Theorems}
\crefname{Lem}{Lemma}{Lemmas}
\crefname{Prop}{Proposition}{Propositions}
\crefname{Cor}{Corollary}{Corollaries}

\author[J. Boulanger]{Julien Boulanger}
\address[Julien Boulanger]{Centro de Modelamiento Matemático (CNRS IRL2807), Universidad de Chile, Santiago, Chile}
\email{jboulanger@cmm.uchile.cl}
\urladdr{https://julien-boulanger.github.io/webpage/}

\author[R. Gutiérrez-Romo]{Rodolfo Gutiérrez-Romo}
\address[Rodolfo Gutiérrez-Romo]{Departamento de Ingeniería Matemática, Facultad de Ciencias Físicas y Matemáticas, Universidad de Chile \& Centro de Modelamiento Matemático (CNRS IRL2807), Universidad de Chile, Santiago, Chile}
\email{g-r@rodol.fo}
\urladdr{http://rodol.fo}

\author[E. Lanneau]{Erwan Lanneau}
\address[Erwan Lanneau]{UMR CNRS 5582, Univ.\ Grenoble Alpes, CNRS, Institut Fourier, F-38000 Grenoble, France}
\email{erwan.lanneau@univ-grenoble-alpes.fr}
\urladdr{https://www-fourier.ujf-grenoble.fr/~lanneau}


\DeclareMathOperator{\lcm}{lcm}

\makeatletter
\g@addto@macro\bfseries{\boldmath}
\makeatother

\usepackage[style=alphabetic,natbib=true,giveninits=true,url=false,maxbibnames=99]{biblatex}
\DeclareFieldFormat{title}{\mkbibquote{#1}}
\AtEveryBibitem{\clearfield{issn}}
\AtEveryCitekey{\clearfield{issn}}

\begin{document}
	
	\begin{abstract}
		Fix $g \geq 2$. Let $\translations(g)$ be the maximal order of the translation group among all genus-$g$ abelian differentials. By work of Schlage-Puchta and Weitze-Schmithüsen, $\translations(g) \leq 4(g - 1)$. They also classify the $g$ attaining this bound. We assume $g$ is outside this class.
		
		We first prove that either $\translations(g) = (2(m + 1) / m) (g - 1)$ for some $m \in \NN \setminus \{0\}$, when regular genus-$g$ origamis exist, or $\translations(g) = 2(g - 1)$, when they do not exist.
		
		In the former case, only some values of $m > 1$ are realizable; $m = 5$ is the smallest. The resulting set of genera, those satisfying $\translations(g) = (12/5)(g - 1)$, contains infinitely long arithmetic progressions. The same holds for any odd prime $m$ congruent to $2$ modulo $3$.
		
		In the latter case, ``many'' strata of the form $\calH(g - 1, g - 1)$, $\calH(2k^q)$ or $\calH(k^{2q})$, where $k \geq 1$ is an integer and $q$ is prime, contain no regular origamis; we derive a complete classification. As an application, we exhibit infinite families of genera $g$ for which $\translations(g) = 2(g - 1)$: $g = p + 1$ for prime $p \geq 5$; $g = p^2 + 1$ for prime, but not Sophie Germain prime, $p$; and $g = pq + 1$, for distinct primes $p, q \geq 5$.
	\end{abstract}
	
	\maketitle
	
	\section{Introduction}
	
	In 1892, Hurwitz \cite{Hurwitz:upper_bound} obtained a celebrated upper bound of $84(g - 1)$ on the maximal number of automorphisms a genus-$g$ compact Riemann surface may possess. Several decades later, between 1968 and 1969, Accola \cite{Accola:lower_bound} and Maclachlan \cite{Maclachlan:lower_bound} independently obtained a lower bound of $8(g + 1)$. Both bounds are known to be sharp in the sense that they are attained for infinitely many $g$.
	
	A finer question is: what is the maximal number of automorphisms of genus-$g$ compact Riemann surfaces? The answer is known for small $g$, and for several infinite families of genera, although it is not known in full generality~\cite{Kiley, Belolipetsky_Jones_p+1, order, May_Zimmerman:automorphisms}.
	
	A related question is counting automorphisms with additional restrictions. This subject has been extensively studied and several versions have been considered, such as counting automorphisms of specific families of compact Riemann surfaces (e.g.\ $p$-gonal surfaces \cite{COSTA_IZQUIERDO,BARTOLINI_P_Gonal}, pseudo-real surfaces \cite{Bujalance_pseudo_real}, or others \cite{Izquierdo_Reyes-Carocca_Rojas}), or those with a prescribed group structure (e.g.\ cyclic or abelian automorphism groups \cite{Wiman_cyclic_automorphism_groups, Harvey_cyclic_automorphism_groups, Maclachlan_abelian_groups, Hidalgo_others_Quasi_Abelian_automorphisms} or with a specific number of elements \cite{Kulkarni:infinite_families,Carocca_Reyes-Carocca, Izquierdo_Reyes-Carocca_Rojas}).
	
	This article focuses on the automorphisms of a compact Riemann surface that also preserve a given holomorphic $1$-form, that is, an \emph{abelian differential}.
	
	A Riemann surface $X$ endowed with a nonzero abelian differential $\omega$ is called a \emph{translation surface} \cite{AM_book, DHV_book}. For this reason, we use the term \emph{translation group} of $(X,\omega)$ for the subgroup of $\Aut(X)$ that preserves $\omega$, that is,
	\[
	\trans(X, \omega) = \{f \in \Aut(X) \st f^*\omega = \omega\}.
	\]
	By slightly abusing notation, we often omit the differential $\omega$ and refer to the pair $(X,\omega)$ simply as $X$ and to its translation group as $\trans(X)$. We will always assume $X$ to be compact.
	
	More precisely, our aim is to study the quantity
	\[
	\translations(g) = \sup \left\{|\trans(X, \omega)| \ \mathrel{\Big|}\ 
	\begin{array}{c}
		\text{$X$ is a genus-$g$ compact Riemann surface} \\
		\text{and $\omega$ is a nonzero abelian differential on $X$}
	\end{array}
	\right\}.
	\]
	The investigation on this quantity was initiated in 2017 by Schlage-Puchta and Weitze-Schmithüsen \cite{Puchta_Schmithusen}. They first show that $\translations(g) \leq 4(g - 1)$ for every $g \geq 2$. In addition, they prove that the pairs $(X,\omega)$ attaining this bound are essentially regular covers, branched over a single point, of pairs $(Y,\eta)$ for some $Y$ of genus $1$. Such pairs $(X,\omega)$ are sometimes referred to as \emph{regular origamis}. More precisely, they prove, using the Riemann--Hurwitz formula, that the upper bound $4(g-1)$ is attained if and only if $(X,\omega)$ is a regular origami and the abelian differential $\omega$ has exactly $2g-2$ zeros of order one. Finally, they completely characterize the $g \geq 2$ such that $\translations(g) = 4(g - 1)$ as those with $g - 1$ divisible by $2$ or $3$.
	
	Our goal is to generalize these results. By refining their application of the Riemann--Hurwitz formula, we show:
	\begin{ThmA}[{\Cref{thm:not_regular}}]\label{thmA:t(g)}
		Let $g \geq 2$. The number $\translations(g)$ is always of the form
		\[
		\translations(g) = \constant(g)(g - 1).
		\]
		where the ``slope'' $\constant(g)$ is either $2$, or has the form $2(m+1)/m$ for some integer $m = \orderm(g) \geq 1$ such that $m \divides 2(g - 1)$, $3 \ndivides m$, and $4 \ndivides m$.
		
		Finally, the case $\constant(g)=2$ arises if and only if no genus-$g$ regular origamis exist.
	\end{ThmA}
	
	In particular, the number $\constant(g)$ belongs to the set
	\[
	\left\{2 < \dotsb < \frac{36}{17} < \frac{15}{7} < \frac{28}{13} < \frac{24}{11} < \frac{11}{5} < \frac{16}{7} < \frac{12}{5} < 3 < 4\right\}.
	\]
	We will see that some of these numbers, such as $3$ and $16/7$, do not actually occur.
	
	\begin{Notation}
		In the case where $\constant(g) = 2$, we will say that $\orderm(g) = \infty$. This notation is justified by the fact that $\lim_{m \to \infty} 2(m + 1)/m = \inf_{m \geq 1} 2(m + 1)/m = 2$.
		
		Moreover, for (possibly infinite) $m \geq 1$, we define the set $\genera(m)$ of genera $g \geq 2$ for which $\orderm(g) = m$, that is,
		\[
		\genera(m) = \{g \geq 2 \st\orderm(g) = m\} = \begin{cases}
			\displaystyle\left\{g \geq 2 \ \mathrel{\Big|}\ \translations(g) = \frac{2(m + 1)}{m}(g - 1)\right\} & \text{if $m < \infty$} \\
			\{g \geq 2 \st\translations(g) = 2(g - 1)\} & \text{if $m = \infty$.}
		\end{cases}
		\]
	\end{Notation}
	Observe that $\genera(1)$ is completely characterized by the work of Schlage-Puchta and Weitze-Schmithüsen. Moreover, \Cref{thmA:t(g)} implies that $\genera(m)$ is empty when $3 \divides m$ or $4 \divides m$. Our aim is to study the sets $\genera(m)$ for $m$ outside of these cases.
	
	We start by determining that $\genera(m)$ is empty for some particular values of $m$.
	
	\begin{ThmA}[\Cref{thm:G(2)_empty,thm:G(2^k-1)_empty}]\label{thm:G_m_empty}
		The sets $\genera(2)$ and $\genera(2^\alpha - 1)$ are empty for each $\alpha \geq 2$.
	\end{ThmA}
	
	As a consequence, the smallest integer $m$, outside of the case $m = 1$ \cite{Puchta_Schmithusen}, such that $\genera(m)$ is nonempty, satisfies $m \geq 5$. We prove that $\genera(5)$ is infinitely large; the proof also works for other prime $m$:
	
	\begin{ThmA}[{\Cref{thm:G_m_general}}] \label{thm:G_m_intro}
		If $m$ is an odd prime and is congruent to $2$ modulo $3$, then $\genera(m)$ contains explicit infinitely long arithmetic progressions. In particular, $\genera(5)$ contains all $g \geq 2$ of the form
		\[
		g = \frac{5k p(p - 1)(p + 1)}{24} + 1
		\]
		where $p$ is prime, $p \bmod{72}$ is $11$, $13$, $59$, or $61$, and $k \equiv \pm 1 \mod{6}$.
	\end{ThmA}
	
	To prove \Cref{thm:G_m_intro}, we construct regular origamis with a translation group of the form $\PSL(2, p) \times \ZZ/k\ZZ$ and show that they achieve the maximal number of translations in genus $g$.
	
	\begin{Rema}
		Numerical experiments suggest that the density of $\genera(5)$ is well-defined and at least $0.001957$.
	\end{Rema}
	
	Recall now that a \emph{Sophie Germain prime} is a prime $p$ such that $2p + 1$ is also prime. Our next result shows that the set $\genera(\infty)$ is also infinitely large:
	
	\begin{ThmA}[{\Cref{thm:no_regular_origamis}}]\label{thm:no_regular_origamis_intro}
		We have that $g \in \genera(\infty)$ for each $g$ of any of the following forms:
		\begin{itemize}
			\item $g = p + 1$, where $p \geq 5$ is prime;
			\item $g = p^2 + 1$, where $p$ is prime, but is not a Sophie Germain prime;
			\item $g = pq + 1$, where $p, q \geq 5$ are distinct primes.
		\end{itemize}
	\end{ThmA}
	
	In particular, \Cref{thm:no_regular_origamis_intro} shows that, for infinitely many $g \geq 2$, no genus-$g$ regular origamis exist.
	
	\begin{Rema}
		A follow-up question is if regular origamis exist in genus $g = pqr + 1$ when $p, q, r \geq 5$ are distinct primes. We know that both cases can arise: they exist if $g = 456 = 5 \cdot 7 \cdot 13 + 1$ (from \Cref{thm:G_m_intro} with $m = 5$, $p = 13$, and $k = 1$), and do not exist if $g = 386 = 5 \cdot 7 \cdot 11 + 1$ (from computer experiments).
	\end{Rema}
	
	The previous theorem suggests focusing on the case of Sophie Germain primes, allowing us to show the following:
	
	\begin{ThmA}[{\Cref{thm:no_regular_origamis,cor:p^2l+1}}]\label{thm:SG_intro}
		If $p \geq 5$ is a Sophie Germain prime, then $p^2 + 1 \in \genera(2p)$. Moreover, $\ell p^2 + 1 \notin \genera(\infty)$ for every integer $\ell \geq 1$.
	\end{ThmA}
	
	If $p$ is a Sophie Germain prime, the first part of the statement shows, in particular, that the set $\genera(2p)$ is nonempty. The second part is equivalent to the existence of regular origamis of genus $g = \ell p^2 + 1$.
	
	\begin{Rema}
		Whenever $p$ is a Sophie Germain prime, \Cref{thm:SG_intro} shows, in particular, that regular origamis exist in genus $g = p^\alpha + 1$ for every $\alpha \geq 2$.
		
		Assume now that $p$ is prime, but not a Sophie Germain prime. \Cref{thm:no_regular_origamis_intro} shows that regular origamis do not exist in genus $g = p^2 + 1$. It is natural to ask if they exist in genus $g = p^\alpha + 1$ for $\alpha > 2$. We know that they do not exist, for example, if $p = 7$ and $\alpha = 3$, that is, when $g = 344 = 7^3 + 1$ (from computer experiments).
	\end{Rema}
	
	The proofs of \Cref{thm:no_regular_origamis_intro,thm:SG_intro} are based on analyzing the existence of regular origamis whose associated abelian differentials have zeros of specific types, that is, belonging to particular \emph{strata}. We focus on the cases $g=p+1$ and $g=pq+1$ primarily because, in these situations, the Euler characteristic $2-2g$ admits a simple prime factorization, and the problem reduces to finding groups with a cyclic subgroup of prime or twice-prime index.
	
	Let us point out that the full automorphism group (and some of its subgroups) of a Riemann surface of genus $g = p+1$ and $g=p^2+1$ has been considered before in the literature \cite{Belolipetsky_Jones_p+1, Izquierdo_Jones_Reyes_p+1,Carocca_Reyes-Carocca}. While the group-theoretic setting differs, the simple factorization of the Euler characteristic similarly facilitates the analysis.
	
	\begin{Rema}
		As a consequence of \Cref{thmA:t(g),thm:G_m_empty,thm:G_m_intro,thm:no_regular_origamis_intro,thm:SG_intro}, an updated list of the possible slopes $\constant(g)$ such that $\translations(g) = \constant(g)(g - 1)$ is:
		\[
		\left\{2 < \dotsb < \frac{52}{25} < \frac{48}{23} < \frac{23}{11} < \frac{40}{19} < \frac{36}{17} < \frac{15}{7} < \frac{28}{13} < \frac{24}{11} < \frac{11}{5} < \frac{12}{5} < 4\right\}.
		\]
		We do not know if the slopes $28/13$, $15/7$, $40/19$, or $52/25$ are realizable.
		
	\end{Rema}
	
	\subsection{Context and motivations}
	This article lies at the intersection between the study of automorphism groups of compact Riemann surfaces and the study of abelian differentials. Whereas the study and classification of the automorphism groups of compact Riemann surfaces is a classical problem that has attracted considerable interest since the late nineteenth century, the case of abelian differentials has cemented its relevance only during the last decades, especially in relation to the study of moduli spaces \cite{AM_book, DHV_book}.
	
	\bigbreak
	
	The subgroup $\trans(X,\omega)$ of $\Aut(X)$ associated with an abelian differential $\omega$ has been studied in several recent works \cite{Puchta_Schmithusen, Hidalgo_automorphisms, Flake_Thevis}. As in the general case, every finite group can be achieved as the translation group of a pair $(X,\omega)$ \cite{Hidalgo_automorphisms}. However, when the genus $g$ of $X$ or the orders of the zeros of $\omega$ are prescribed, this group had not, to the best of our knowledge, been investigated outside two cases: the ``upper bound'' case $\translations(g) = 4(g-1)$ \cite{Puchta_Schmithusen}, and the case where it is a $p$-group \cite{Flake_Thevis}.
	
	\subsubsection*{Geometric interpretation} 
	This work was originally motivated by the connection between abelian differentials and flat geometry. As previously mentioned, a Riemann surface $X$ with a nonzero abelian differential $\omega$ is also called a \emph{translation surface}. Translation surfaces admit other equivalent definitions, which also provide equivalent definitions of the translation group. A more combinatorial definition is a collection of polygons on the plane with side identifications by translations up to scissors congruences. Equivalently, it is a genus-$g$ topological surface $S$ endowed with a \emph{translation atlas}, that is, an atlas whose transition functions are translations, except at finitely many points called singularities. This atlas allows us to define the total area of $X$. Since we assume $X$ to be compact, this area is finite.
	
	Using the translation atlas, we may also define the group of \emph{affine homeomorphisms} $\Aff(X,\omega)$ of a translation surface as the subgroup of $\Homeo^+(S)$ with constant derivative in the atlas. This group is well-defined since any matrix remains constant when conjugated by a translation.
	By taking the derivative of an affine homeomorphism, we obtain the \emph{derivative map} $D \colon \Aff(X,\omega) \to \SL(2, \RR)$. The image of this map is known as the \emph{Veech group} $\SL(X,\omega)$ of $(X,\omega)$ and records all possible matrices that can be lifted to affine homeomorphisms. The kernel of this map is exactly the translation group $\trans(X,\omega)$.
	
	\subsubsection*{Regular origamis} One of the simplest examples of a translation surface is the unit square torus $\TT = \RR^2/\ZZ^2$, equipped with the $1$-form $\mathrm{d}z$.
	By considering covers of $\TT$ branched over a single point (given by the points of integer coordinates), we obtain an \emph{origami} or \emph{square-tiled surface}. The differential also lifts to the covering surface.
	
	A particular case of origamis of special interest is \emph{regular origamis} (also known as \emph{normal origamis}): those for which the cover to the unit torus is normal. Regular origamis can also be defined in terms of their translation group. Indeed, given a finite group $G$ generated by two elements $x$ and $y$, we can define a regular origami by labeling unit squares with the elements of $G$, and declaring that rightward gluings are given by multiplication by $x$, and upward gluings, by $y$. Then, the translation group of the resulting origami is isomorphic to $G$.
	
	It turns out regular origamis constitute the translation surfaces with the largest translation group. Namely, a translation surface is a regular origami if and only if its translation group has more than $2(g-1)$ elements (see \Cref{lem:origami_bounds}). Furthermore, if $(X, \omega)$ is a regular origami, the order $\trans(X, \omega)$ only depends on the order $m$ of the zeros of $\omega$. Indeed, it equals $(2(m + 1)/m)(g - 1)$. This number can be computed from the generators $x, y \in G$ as $m = \ord([x, y]) - 1$. See \Cref{sec:buiding_regular_origamis}.
	
	\subsubsection*{Strata of abelian differentials.}
	In geometric terms, the zeros of the abelian differential $\omega$ correspond to the singularities of the translation surface $(X, \omega)$.
	In fact, the (moduli) space of translation surfaces is partitioned into \emph{strata} prescribing the orders of the zeros of $\omega$.
	
	\begin{Notation}
		As standard in the theory, we will denote the set of genus-$g$ translation surfaces whose abelian differential has $s_i$ zeros of order $k_i$, for $i \in \{1, \dotsc, \ell\}$, by $\calH_g(k_1^{s_1}, \dots, k_\ell^{s_\ell})$. As such, a superscript in this notation will always mean a \emph{multiplicity} (and never an exponent). We refer to such a set as a \emph{stratum}.
		
		The Riemann--Roch Theorem relates the singularity data and the genus:
		\begin{equation}\label{eq:Riemann_Roch}
			2g-2 = \sum_{i=1}^\ell s_ik_i.
		\end{equation}
		Hence, we often omit the subscript $g$.
	\end{Notation}

	The work of Schlage-Puchta and Weitze-Schmithüsen \cite{Puchta_Schmithusen} focuses on regular origami in the stratum $\calH(1^{2g-2})$. To obtain the first part of \Cref{thm:G_m_empty}, and \Cref{thm:no_regular_origamis_intro,thm:SG_intro}, we study the existence of regular origamis in several other strata. We obtain:
	\begin{ThmA}\label{thm:all_strata}
		Let $k, \ell \geq 1$ be integers. Then:
		\begin{itemize}
			\item $\calH(k^\ell)$, for even $k, \ell$, contains regular origamis (\Cref{sec:examples_origamis} Example \eqref{i:example_2k^2l});
			\item $\calH(2^\ell)$, for odd $\ell$, contains regular origamis if and only if $\ell$ is divisible by $9$ (\Cref{thm:H(2^g-1)});
			\item $\calH(k^2)$, for odd $k$, contains no regular origamis (\Cref{lem:H(g-1 g-1)});
			\item $\calH(k^4)$ contains regular origamis (\Cref{sec:examples_origamis} Example \eqref{i:example_k^4});
			\item $\calH(k^6)$, when $k \equiv 1 \mod{4}$, contains regular origamis if and only if every prime factor of $(k + 1) / 2$ is congruent to $1$ modulo $3$ (\Cref{thm:stratum_H(k^2q)}); and
			\item $\calH(k^6)$, when $k \equiv 3 \mod{4}$, contains regular origamis if and only if every prime factor of $(k + 1) / 4$ is congruent to $1$ modulo $3$ (\Cref{thm:stratum_H(k^2q)});
		\end{itemize}
		Furthermore, if $q$ is an odd prime, then:
		\begin{itemize}
			\item $\calH(k^{2q})$, for odd $k$ and $q > 3$, contains no regular origamis (\Cref{thm:stratum_H(k^2q)}); and
			\item  $\calH(2k^q)$ contains regular origamis if and only if every prime factor of $2k + 1$ is congruent to $1$ modulo $q$ (\Cref{thm:stratum_H(2p^q)}).
		\end{itemize}
	\end{ThmA}
	
	\begin{Rema}\label{rema:strata}
		From \Cref{eq:Riemann_Roch}, translation surfaces in $\calH(k^\ell)$ have genus $g$ satisfying $2g - 2 = k\ell$. In particular, $k$ and $\ell$ cannot both be odd.
		
		In the case of $\ell = 2$, we have $\calH(k^2) = \calH(g - 1, g - 1)$ for genus $g = k + 1$. The previous theorem states that this stratum contains regular origamis if and only if $g$ is odd. Similarly, surfaces in $\calH(2^\ell)$ have genus $g = \ell + 1$. Thus, this stratum contains regular origamis if and only if $g - 1$ is divisible by $2$ or $9$.
	\end{Rema}
	
	This result is somewhat complementary to the work of Flake and Thevis \cite{Flake_Thevis}, which investigates the strata in which a regular origamis whose translation group is a $p$-group may occur. Together with the classification of $\genera(1)$ \cite{Puchta_Schmithusen}, their work shows that $1 < \orderm(g) \leq p^\alpha - 1$ for every $g$ of the form $g = p^\beta(p^\alpha - 1)/2 + 1$, where $p > 3$ is prime, $\alpha \geq 1$ and $\beta \geq \alpha + 1$. In contrast, \Cref{thm:all_strata} does not make assumptions about the form of the group, but only deals with specific strata.
	
	It is known that a regular origami constructed from a group $G$ and two generators $x,y$ belongs to the stratum $\calH(k^\ell)$ if and only if the cyclic subgroup $H = \langle [x,y] \rangle$ has order $k+1$ and index $\ell$ in $G$ \cite{Puchta_Schmithusen, Flake_Thevis} (see \Cref{sec:buiding_regular_origamis}). As a consequence, \Cref{thm:all_strata} reduces to a classification problem for groups of order $(k+1)\ell$ generated by two elements whose commutator has order $k+1$. In the case where $k+1$ and $\ell$ have simple prime factorizations, we are able to provide a full classification.
	
	Furthermore, we can further classify the translation groups of regular origamis in the strata $\calH(k^6)$ when $k$ is odd, and $\calH(2k^q)$. Indeed, the former case only admits groups of the form $(\ZZ/\lambda\ZZ \times \ZZ/2\ZZ \times \ZZ/2\ZZ) \rtimes \ZZ/3\ZZ$, or $(\ZZ/\lambda\ZZ \times \Q_8) \rtimes \ZZ/3\ZZ$; the latter case only admits groups of the form $\ZZ/(2k + 1)\ZZ \rtimes \ZZ/q\ZZ$. See \Cref{thm:stratum_H(2p^q),thm:stratum_H(k^2q)} for more details.
	
	\subsection*{Acknowledgments} The first author was supported by Centro de Modelamiento Matemático (CMM) BASAL fund FB210005 for center of excellence from ANID-Chile.
	
	The second author is grateful to Ferrán Valdez for introducing them to the original article by Schlage-Puchta and Weitze-Schmithüsen, to Gabriela Weitze-Schmithüsen for interesting discussions, and is supported by Centro de Modelamiento Matemático (CMM) BASAL fund FB210005 for center of excellence from ANID-Chile, and the FONDECYT Regular N°1221934 and N°1250798 grants from ANID-Chile.
	
	The third author has been partially supported by the LabEx PERSYVAL-Lab (ANR-11-LABX-0025-01) funded by the French program Investissement d’avenir and by the CNRS. The third author would like to thank the Center for Mathematical Modeling (CMM) at Universidad de Chile for excellent working conditions.
	
	We are also grateful to Milagros Izquierdo and Sebastián Reyes-Carocca for pointing us to interesting references on automorphism groups of Riemann surfaces, and to David Grimm for the beautiful argument in the proof of \Cref{lem:sum_of_two_squares}.
	
	\section{Background and preliminaries}
	
	In this section, we will first provide the necessary context in group theory. We will also state and prove a series of simple lemmas that will be useful later.
	
	\subsection{Basic facts about groups}
	
	Throughout the article, we will mainly deal with finite groups. Given a finite group $G$, its \emph{order} is its number of elements, which will be denoted by $|G|$. If $H \leq G$ is a subgroup of $G$, its \emph{index} is the number of cosets of $H$ inside $G$, equals $|G|/|H|$, and is denoted by $(G : H)$ (other common notations include $[G : H]$ and $|G : H|$.
	
	We start with some well-known facts about subgroups.
	
	\begin{Lem} \label{lem:order_HK}
		Let $G$ be a finite group and let $H, K \leq G$ be subgroups of $G$. We have that
		\[
		|H K| = \frac{|H| |K|}{|H \cap K|}.
		\]
	\end{Lem}
	
	\begin{proof}
		Consider the group $H \times K$ acting on the set $H K$ via $(h, k) x = h x k^{-1}$. The action is transitive and the stabilizer of $1 \in H K$ is isomorphic to $H \cap K$. Hence, by the orbit-stabilizer theorem, we get
		\[
		|H K|\cdot |H \cap K| = |H \times K| = |H| |K|.
		\]
		Solving for $|H K|$ yields the desired result.
	\end{proof}
	
	\begin{Lem}\label{lem:index_two_normal}
		Let $G$ be a group and let $H \leq G$ be an index-two subgroup. Then, $H$ is normal in $G$.
	\end{Lem}
	\begin{Lem} \label{lem:index_two_order_two}
		Let $G$ be a group and let $H \leq G$ be an index-two subgroup. If $x, y \in G \setminus H$, then $xy \in H$. In particular, $x^2 \in H$ for every $x \in G$.
	\end{Lem}
	\begin{proof}
		The group $H$ is normal by the previous lemma. Thus, if $x, y \in G \setminus H$, we have that $x$ and $y$ project to the only nontrivial element of $G / H \simeq \ZZ/2\ZZ$. Hence, $x y H = H$, so $x y \in H$.
		
		Finally, if $x \in G \setminus H$, we get that $x^2 \in H$. If $x \in H$, we also have $x^2 \in H$.
	\end{proof}
	
	We recall that a subgroup $H \leq G$ is \emph{characteristic} if it is preserved (setwise) by every automorphism of $G$, and continue with a simple fact stating that a subgroup of a cyclic group is cyclic and uniquely determined by its order:
	\begin{Thm}\label{thm:fundamental_cyclic_groups}
		Let $G$ be a finite cyclic group. We have that every subgroup of $G$ is cyclic. Moreover, there exists a unique such subgroup of order $k$ for every divisor $k$ of $|G|$. In particular, every subgroup of $G$ is characteristic.
	\end{Thm}
	
	\begin{Lem}\label{lem:characteristic_in_normal_is_normal}
		Let $G$ be a group, and assume that $H$ is a normal subgroup of $G$. If $K$ is a characteristic subgroup of $H$, then $K$ is normal in $G$.
	\end{Lem}
	\begin{proof}
		Consider the action $\varphi \colon G \to \Aut(H)$ given by conjugation. This morphism is well-defined since $H$ is normal in $G$. Since $K$ is characteristic in $H$, we deduce that $\varphi(g)$ stabilizes $K$ for every $g \in G$. Thus, $K \triangleleft G$.
	\end{proof}
	
	If $g, h \in G$, we denote their \emph{commutator} by $[g, h] = g h g^{-1} h ^{-1}$. Recall that the \emph{commutator subgroup} of $G$ (also known as the \emph{derived subgroup} of $G$) is the group generated by its commutators. We denote it by $G'$ or $[G,G]$. The commutator subgroup is characteristic, and $G / H$ is abelian whenever $G' \leq H \triangleleft G$ (particularly when $H = G'$).
	
	Now, it is possible to iterate the derivation process and consider the \emph{derived series} associated with $G$:
	\[
	G^{(0)} = G, \quad G^{(1)} = [G^{(0)},G^{(0)}], \quad G^{(2)} = [G^{(1)},G^{(1)}], \quad \dots
	\]
	
	A group is called \emph{solvable} if its derived series eventually reaches the trivial group. The celebrated results of Burnside and Feit--Thomson give criteria on the order of a group for it to be solvable:
	
	\begin{Thm}[{Burnside's theorem \cites{Burnside:p^a-q^b}[Theorem 7.8]{Isaacs_Gp_theory}}] \label{thm:Burnside-p^a_q^b}
		If $p$, $q$ are prime and $\alpha, \beta$ are nonnegative integers, then every group of order $p^{\alpha} q^{\beta}$ is solvable.
	\end{Thm}
	
	\begin{Thm}[{Feit--Thompson \cite{Feit-Thompson1,Feit-Thompson2}}] \label{thm:feit-thompson}
		Any finite group of odd order is solvable.
	\end{Thm}
	
	Roughly speaking, solvable groups are those that can be split into abelian blocks. Examples of nonsolvable groups include (nontrivial) \emph{perfect groups}: groups $G$ with $G'=G$. These results show that perfect groups can only exist for some orders.
	
	Recall that the \emph{center} of $G$, denoted $\Z(G)$, is the subgroup of those elements of $G$ commuting with every other element of $G$. We include two facts about the center of a perfect group for later use:
	\begin{Lem}\label{lm:action}
		Assume that $H$ is a cyclic normal subgroup of a finite group $G$. Then, $H \leq \Z(G')$. In particular, if $G$ is perfect, then $H \leq \Z(G)$.
	\end{Lem}
	\begin{proof}
		Since $H$ is normal, the action by $G$ on $H$ by conjugation is well-defined, and induces a homomorphism $G \to \Aut(H)$. Since $H$ is cyclic, the group $\Aut(H)$ is abelian and, thus, this homomorphism factors through the abelianization $G/G'$. 
		In other words, $H$ is contained in its kernel, namely $H \leq \Z(G')$.
	\end{proof}
	
	\begin{Lem}[{Grün's Lemma \cite[p.\ 61]{Rose_groups}}]\label{lem:Grun's}
		If $G$ is a perfect group, then $G/\Z(G)$ has a trivial center.
	\end{Lem}
	
	\subsection{Semidirect products and their commutator subgroup.}
	Recall that a \emph{semidirect product} between groups $N$ and $H$ via a homomorphism $\varphi \colon H \to \Aut(N)$ is the group $G$ of with underlying set $N \times H$ and the operation
	\[
	(n_1, h_1) \cdot (n_2, h_2) = (n_1 \varphi(h_1)(n_2), h_1 h_2).
	\]
	We denote it by $G = N \rtimes_{\varphi} H$. We will sometimes omit the map $\varphi$.
	
	We will encounter commutator subgroups of semidirect products. The following lemma follows directly from the definitions:
	\begin{Lem}\label{lem:derived_sbgp_semi_direct}
		The commutator subgroup of a semidirect product $G = N \rtimes H$ is a subgroup of $N \rtimes H'$, where $H'$ is the commutator subgroup of $H$.
	\end{Lem}
	
	A particular case is a semidirect product of cyclic groups. If $k, \ell$ are integers and $d^\ell \equiv 1 \mod{k}$, we can define $\varphi_d \colon \ZZ/\ell\ZZ \to \Aut(\ZZ/k\ZZ)$ by declaring that $\varphi_d(1)$ maps $1$ to $d$, and extending by cyclicity. We denote $\ZZ/k\ZZ \rtimes_{\varphi_d} \ZZ/\ell\ZZ$ simply as $\ZZ/k\ZZ \rtimes_d \ZZ/\ell\ZZ$. More explicitly:
	\[
	(a_1, b_1) \cdot (a_2, b_2) = (a_1 + d^{b_1} a_2, b_1 + b_2).
	\]
	
	\begin{Lem}
		In the context above, $\ZZ/k\ZZ\rtimes_d \ZZ/\ell\ZZ$ has a presentation:
		\[
		\langle x,y \st x^{k} = y^\ell = 1 \text{ and } [y, x] = x^{d-1} \rangle.
		\]
	\end{Lem}
	
	\begin{proof}
		Take $x = (1, 0)$ and $y = (0, 1)$.
	\end{proof}
	
	Semidirect products of cyclic groups fall under the more general class of \emph{metacyclic groups}. These are defined as extensions of cyclic groups by cyclic groups, meaning that they obey a short exact sequence
	\[
	1 \to \ZZ/k\ZZ \to G \to \ZZ/\ell\ZZ \to 1,
	\]
	and have a presentation
	\[
	G = \langle x, y \st  x^k =1,  y^\ell= x^r \text{ and } [y, x] = x^{d - 1} \rangle,
	\]
	where, as before, $d^\ell = 1 \mod{k}$ and, moreover, $k \divides r(d - 1)$ and $r \divides k$. In this context, a metacyclic group is a semidirect product of cyclic groups if and only if $r = k$.
	Every element of a metacyclic group admits a normal form: it can be written uniquely as $y^\alpha x^\beta$ for $0 \leq \alpha < \ell$ and $0 \leq \beta < k$. Indeed, the group $K = \langle x \rangle$ is normal, cyclic of order $k$, and $G/K$ has order $\ell$. Indeed, the cosets $y^\alpha (G/K)$ are distinct for distinct $\alpha$.
	
	\medbreak
	
	The following lemma is well-known:
	\begin{Lem}\label{lem:commutator_metacyclic}
		Using the presentation above, the commutator subgroup of a metacyclic group $G$ is $G' = \langle x^{d-1} \rangle$. In the particular case of $G = \ZZ/k\ZZ \rtimes_d \ZZ/\ell\ZZ$, we have $G' \simeq \ZZ/t\ZZ \leq \ZZ/k\ZZ$, where $t = k / \gcd(d - 1, k)$.
	\end{Lem}
	\begin{proof}
		Let $H = \langle x^{d-1} \rangle$. Since $x^{d-1} = [y, x]$, we have $H \leq G'$. We will show that $G' \leq H$.
		
		We start by showing that, if $\alpha,\beta \in \ZZ$, then:
		\[
		y^\alpha x^\beta y^{-\alpha} = x^{\beta d^\alpha}.
		\]
		First, from the relation $[y, x] = x^{d-1}$ we get $y x y^{-1} = x^d$, so
		\[
		y^\alpha x y^{-\alpha} = y^{\alpha-1}x^dy^{-\alpha+1} = \left(y^{\alpha-1}xy^{-\alpha+1}\right)^d.
		\]
		Inductively,
		\[
		y^\alpha x y^{-\alpha} = x^{d^\alpha}
		\]
		and, therefore,
		\[
		y^\alpha x^\beta y^{-\alpha} = \left(y^\alpha x y^{-\alpha}\right)^\beta = x^{\beta d^\alpha}.
		\]
		Now, we have that
		\[
		[y^\alpha, x^\beta] = (y^\alpha x^\beta y^{-\alpha}) x^{-\beta} = x^{\beta d^\alpha} x^{-\beta} = x^{\beta(d^\alpha - 1)} = (x^{d^\alpha - 1})^\beta.
		\]
		Moreover,
		\[
		x^{d^\alpha - 1} = (x^{d - 1})^{1 + d + d^2 + \dotsb + d^{\alpha-1}} \in H,
		\]
		so $[y^\alpha, x^\beta] \in H$.
		
		Finally, if $\alpha_1, \beta_1, \alpha_2, \beta_2 \in \ZZ$, we have:
		
		\begin{align*}
			[y^{\alpha_1} x^{\beta_1}, y^{\alpha_2} x^{\beta_2}] &= y^{\alpha_1} x^{\beta_1} y^{\alpha_2} x^{\beta_2} x^{-\beta_1} y^{-\alpha_1} x^{-\beta_2} y^{-\alpha_2} \\
			&= (y^{\alpha_1} x^{\beta_1} y^{-\alpha_1}) (y^{\alpha_1 + \alpha_2} x^{\beta_2 - \beta_1} y^{-(\alpha_1 + \alpha_2)}) (y^{\alpha_2} x^{-\beta_2} y^{-\alpha_2}) \\
			&= x^{\beta_1 d^{\alpha_1}} x^{(\beta_2 - \beta_1) d^{\alpha_1 + \alpha_2}} x^{-\beta_2 d^{\alpha_2}} \\
			&= x^{-(\beta_1 d^{\alpha_1}) (d^{\alpha_2} - 1)} x^{(\beta_2 d^{\alpha_2}) (d^{\alpha_1} - 1)} \\
			&= \left(x^{d-1}\right)^{(-\beta_1 d^{\alpha_1})(1+d+\cdots+d^{\alpha_2-1})+(\beta_2 d^{\alpha_2})(1+d+\cdots+d^{\alpha_1-1})} \in H.
		\end{align*}
		
		In the particular case where $G = \ZZ/k\ZZ \rtimes \ZZ/\ell\ZZ$, we get that $G' = \langle (d - 1, 0) \rangle$. The order of $(d - 1, 0)$ is exactly $t = k / \gcd(d - 1, k)$, so $G' \simeq \ZZ/t\ZZ$.
	\end{proof}
	
	A fundamental tool in the theory of finite groups is the following result of Schur and Zassenhaus, which in some cases exhibits a group as a semidirect product.
	\begin{Thm}[{Schur--Zassenhaus \cite[Section 3B]{Isaacs_Gp_theory}}]\label{thm:Schur_Zassenhaus}
		Let $G$ be a finite group. Let $H \triangleleft G$ be a normal subgroup of $G$ whose order is coprime to its index in $G$. Then, $G$ can be written as a semidirect product $G \simeq H \rtimes G/H$.
	\end{Thm}
	
	In view of this result, it is natural to look for normal subgroups of prime order or normal subgroups whose index is a maximal power of a prime. This is the content of the next subsection.
	
	\subsection{Basic facts about \texorpdfstring{$p$}{p}-groups} Given a prime $p$, a finite group $G$ is a \emph{$p$-group} if its order is $p^\alpha$ for some $\alpha \in \NN$. We will use several facts about $p$-groups.
	
	\begin{Lem}[{\cite[Corollary 1.24]{Isaacs_Gp_theory}}] \label{lem:p-group_normal_subgroup}
		If $G$ is a $p$-group of order $p^\alpha$, then $G$ contains a normal subgroup of order $p^\beta$ for each $0 \leq \beta \leq \alpha$.
	\end{Lem}
	This allows us to prove:
	
	\begin{Cor}\label{cor:p-group_index_commutator_subgroup}
		If $G$ is a $p$-group of order at least $p^2$, then $(G : G') \geq p^2$.
	\end{Cor}
	
	\begin{proof}
		Assume that the order of $G$ is $p^\alpha$ for $\alpha \geq 2$. By \Cref{lem:p-group_normal_subgroup}, there exists a normal subgroup $H$ of $G$ of order $p^{\alpha - 2}$.
		In particular $(G : H) = p^2$. The group $G/H$ has order $p^2$, so it is abelian. Since $G'$ is the smallest normal subgroup of $G$ such that $G/G'$ is abelian, we get $G' \leq H$. Consequently, $(G : G') \geq (G : H) = p^2$.
	\end{proof}
	
	The study of $p$-groups is central to understanding the structure of finite groups. A cornerstone result is Sylow's theorems. They state that every finite group has subgroups with a maximal prime-power order and derive some of their properties. Given a finite group $G$ and a prime number $p$, write $|G| = p^\alpha k$ for $k$ coprime to $p$. Then, a subgroup of $G$ of order $p^\alpha$ is called a \emph{Sylow $p$-subgroup} of $G$. We have:
	
	\begin{Thm}\cite[Theorems 1.7,1.12, 1.17]{Isaacs_Gp_theory}
		\label{thm:Sylow}
		Let $G$ be a finite group and $p$ be a prime number. Then,
		\begin{itemize}
			\item there exists at least one Sylow $p$-subgroup of $G$;
			\item all such groups are conjugate; and
			\item the number $n_p$ of these groups satisfies $n_p \equiv 1 \mod p$.
		\end{itemize}
	\end{Thm}
	In fact, Hall generalized this result to collections of prime numbers in the case where $G$ is solvable.
	More precisely, given a collection $\pi$ of prime numbers dividing $|G|$, a \emph{$\pi$-Hall subgroup} of $G$ is a subgroup $H$ whose order is a multiple of every prime in $\pi$ and whose index is coprime to every prime in $\pi$. Hall~\cite{Hall_groups} showed that every solvable group contains a $\pi$-Hall subgroup, namely:
	
	\begin{Thm}\cite[Theorem 3.13]{Isaacs_Gp_theory}\label{thm:Hall}
		Suppose $G$ is a finite solvable group and let $\pi$ be a collection of primes dividing $|G|$. Then:
		\begin{itemize}
			\item there exists a $\pi$-Hall subgroup of $G$; and
			\item all such groups are conjugate.
		\end{itemize}
	\end{Thm}
	Furthermore, the number of Hall subgroups is of the form 
	\[ 1 + \sum_{p \in \pi} a_p p\]
	for some integers $a_p \geq 0$ \cite[Lemma 15]{Puchta_Schmithusen}.
	
	\bigbreak
	
	As previously mentioned, understanding the Sylow or Hall subgroups of a group provides deeper insight into its structure. Many of the groups we consider will contain, for any prime $p$ dividing its order, a cyclic subgroup of index $p$ inside any of their Sylow $p$-subgroups. The $p$-groups containing a cyclic subgroup of index $p$ were classified by Burnside. Here we will state this result for $p=2$.
	
	\begin{Thm}[{\cites[IV.4]{Brown}[Proposition 10.1]{Cziszter_cyclic_index_two}[\S 5.3.4]{Book_Robinson_groups}}]\label{thm:Burnside}
		The only finite $2$-groups containing a cyclic subgroup of index two are:
		\begin{enumerate}
			\item $\ZZ/2^\alpha \ZZ$ for $\alpha \geq 1$;
			\item $\ZZ/2^{\alpha-1} \times \ZZ/p\ZZ$ for $\alpha \geq 2$;
			\item $\M_{2^\alpha} =\ZZ/{2^{\alpha-1}\ZZ} \rtimes_d \ZZ/2\ZZ$, where $d=2^{\alpha-2}+1$, for $\alpha \geq 3$.
			
			\item The dihedral group $\D_{2^{\alpha}} = \ZZ/2^{\alpha-1}\ZZ \rtimes_{-1} \ZZ/2\ZZ$, for $\alpha \geq 3$;
			\item $\SD_{2^{\alpha}} = \ZZ/{2^{\alpha-1}\ZZ} \rtimes_d \ZZ/2\ZZ$, where $d=2^{\alpha-2}-1$, for $\alpha \geq 4$; and
			\item The dicyclic group $\Dic_{2^\alpha}$, for $\alpha \geq 3$.
		\end{enumerate}
	\end{Thm}
	The dicyclic group $\Dic_{2^\alpha}$ is not a semidirect product of cyclic groups, but it is a metacyclic group. Concretely, it has a presentation:
	\[
	\Dic_{2^\alpha} = \langle x,y \st x^{2^{\alpha-1}}=1, y^2 = x^{2^{\alpha-2}} \text{ and } [y, x] = x^{-2} \rangle.
	\]
	In particular, its commutator subgroup is isomorphic to $\ZZ/2^{\alpha-2}\ZZ$. 
	
	We now provide useful results about the groups in \Cref{thm:Burnside}:

	\begin{Prop} \label{lem:Burnside_Aut}
		Assume that $G$ is a finite $2$-group containing a cyclic subgroup of index two. Then, $\Aut(G)$ is a $2$-group unless 
		\begin{itemize}
			\item $G = \ZZ/2\ZZ \times \ZZ/2\ZZ$, for which $\Aut(G) \simeq \Sym_3$ has $6$ elements; or
			\item $G = \Dic_8 = \Q_8$, for which $\Aut(G) \simeq \Sym_4$ has $24$ elements.
		\end{itemize}
	\end{Prop}
	\begin{proof}
		We check each case in \Cref{thm:Burnside}:
		\begin{itemize}
			\item The order of $\Aut(\ZZ/2^\alpha\ZZ)$ is $\varphi(2^\alpha) = 2^{\alpha-1}$ \cite[\S 1.5.5]{Book_Robinson_groups}.
			\item The order of $\Aut(\ZZ/2^{\alpha-1}\ZZ \times \ZZ/2\ZZ)$ is $2^\alpha$ except for $\alpha=2$, for which it is $6$ \cite{ShabaniAttar_Automorphisms}.
			\item The order of $\Aut(\D_{2^{\alpha}}) \simeq \Aut(\Dic_{2^{\alpha}})$ is $2^{2\alpha-1}$, except for $\alpha=3$ for which we have $\Aut(\Q_8) \simeq \Sym_4$ \cites{Walls_automorphism_gp_dihedral}[Exercise 5.3.4]{Book_Robinson_groups}.
			\item The order of $\Aut(\SD_{2^\alpha})$ is $\varphi(2^{\alpha-1})\cdot 2^{\alpha-2} = 2^{2\alpha-4}$ \cite[\S 2.3.2]{Thesis_Marin_Montilla}.
			\item The order of $\Aut(\M_{2^\alpha})$ is $2^\alpha$ \cite{ShabaniAttar_Automorphisms}. \qedhere
		\end{itemize}
	\end{proof}
	
	This specificity for the groups $\ZZ/2\ZZ \times \ZZ/2\ZZ$ and $\Q_8$ also imply:
	\begin{Prop}\label{lem:Burnside_characteristic}
		Assume that $G$ is a finite $2$-group containing a cyclic subgroup of index two. Then, $G$ contains a characteristic subgroup of index two, unless $G \simeq \ZZ/2\ZZ \times \ZZ/2\ZZ$ or $G \simeq \Q_8$.
	\end{Prop}
	\begin{Rema}
		The groups $\ZZ/2\ZZ \times \ZZ/2\ZZ$ and $\Q_8$ both contain three (cyclic) subgroups of index two, which are permuted by the automorphisms of order $3$. In particular, none of these subgroups is characteristic.
	\end{Rema}
	\begin{proof}[Proof of \Cref{lem:Burnside_characteristic}]
		We first use that the number of subgroups of index two of a finite group is given by $n = (G:G^2)-1$ \cite{Nganou_sgps_index_2}. By hypothesis, we know that $n > 0$. In particular, since $G$ is a $2$-group and $(G:G^2)$ divides $|G|$, we deduce that $n$ is odd.
		
		Now, $\Aut(G)$ acts on the set of subgroups of index two. By hypothesis and \Cref{lem:Burnside_Aut}, we know that $\Aut(G)$ is a $2$-group, so the size of each orbit of this action is a power of two by the orbit-stabilizer theorem. As the sum of all these sizes is $n$, which is odd, we deduce that there exists an orbit $\{H\}$ of size one. In other words, the group $H \leq G$ is characteristic.
	\end{proof}
	
	We end this section with a discussion on the existence of \emph{normal $p'$-Hall subgroups}.
	
	\subsection{Frobenius \texorpdfstring{$p$}{p}-complement theorem.}
	From the Schur--Zassenhaus theorem, a group $G$ can be split into a semidirect product if one finds a\emph{normal} Hall subgroup. In the specific case where $\pi$ is the collection of primes dividing $G$, except for the prime $p$, we write $\pi = p'$. A normal $p'$-Hall subgroup is called a \emph{normal $p$-complement}.
	
	As usual, we will denote by $\N_G(X) = \{g \in G \st g X g^{-1} = X\}$ the \emph{normalizer} of the subgroup $X$ in $G$ and by $\C_G(X) = \{g \in G \st\text{$g x = x g$ for every $x \in X$}\}$ the \emph{centralizer} of $X$ in $G$. We now state the Frobenius $p$-complement theorem, which provides the existence of a normal $p$-complement under certain conditions. 
	
	\begin{Thm}[{Frobenius \unboldmath{$p$}-complement theorem \cite[Theorem 5.26]{Isaacs_Gp_theory}}]\label{thm:Frobenius_p_complement}
		Let $p$ be a prime number, and let $G$ be a finite group. The following are equivalent:
		\begin{enumerate}
			\item $G$ has a normal $p$-complement; and
			\item for every $p$-subgroup $X \leq G$, the group $\N_G(X)/\C_G(X)$ is a $p$-group.
		\end{enumerate}
	\end{Thm}
	\begin{Rema}
		The original theorem contains a third equivalent statement, but we will not use it.
	\end{Rema}
	
	\begin{Rema}\label{rk:normalizer_centralizer_automorphism}
		Recall, for later use, that the action of $\N_G(X)$ on $X$ by conjugation induces a monomorphism $\N_G(X)/\C_G(X) \hookrightarrow \Aut(X)$.
	\end{Rema}
	
	We continue with several facts about normal $p$-complements.
	\begin{Lem} \label{lem:p_Sylow_order_p_Hall}
		Let $p$ be a prime number. Let $G$ be a finite group admitting a normal $p$-complement $N$. If $p^2 \ndivides (G : G')$, then $G \simeq N \rtimes \ZZ/p\ZZ$.
	\end{Lem}
	
	\begin{proof}
		By the Schur--Zassenhaus theorem (\Cref{thm:Schur_Zassenhaus}), we have
		\[
		G \simeq N \rtimes L,
		\]
		where $L$ is a Sylow $p$-subgroup of $G$. If $|L| \geq p^2$, by \Cref{lem:derived_sbgp_semi_direct} and \Cref{cor:p-group_index_commutator_subgroup} we deduce that $p^2 \divides (G : N \rtimes L')$, so $p^2 \divides (G : G')$. This contradicts the hypothesis.
		
		Therefore, $L \simeq \ZZ/p\ZZ$ and
		\[
		G \simeq N \rtimes \ZZ/p\ZZ. \qedhere
		\]
	\end{proof}
	
	\begin{Cor} \label{cor:p_Sylow_order_p_cyclic}
		Let $p$ be a prime number. Let $G$ be a finite group containing a cyclic normal subgroup $H$ of index $p$. If $H \leq G'$, then $G \simeq \ZZ/\ell \ZZ \rtimes \ZZ/p\ZZ$, where $\ell = |H|$. Moreover, $p \ndivides \ell$.
	\end{Cor}
	
	\begin{proof}
		Write $|H| = k p^\alpha$ with $p \ndivides k$ and $\alpha \geq 0$. Since $H$ is cyclic, we have that:
		\[
		H \simeq \ZZ/k\ZZ \times \ZZ/p^\alpha\ZZ.
		\]
		Moreover, $\ZZ/k\ZZ \leq H$ is characteristic in $H$ by \Cref{thm:fundamental_cyclic_groups} and, since $H$ is normal in $G$, we deduce that $\ZZ/k\ZZ$ is normal in $G$ by \Cref{lem:characteristic_in_normal_is_normal}. Thus, $\ZZ/k\ZZ$ is a normal $p'$-Hall subgroup of $G$. Since $H \leq G'$, we have $p = (G : H) \geq (G : G')$, so \Cref{lem:p_Sylow_order_p_Hall} shows that 
		\[
		G \simeq \ZZ/k\ZZ \rtimes \ZZ/p\ZZ.
		\]
		We deduce that $\alpha = 0$. Taking $\ell = k$, we obtain the desired conclusion.
	\end{proof}
	
	We will also use the following result:
	
	\begin{Lem} \label{lem:cyclic_index_two}
		Let $G$ be a finite group containing a cyclic $2$-subgroup $H$ such that $(G : H)$ is even, but not divisible by $4$. Then, every nontrivial $2$-subgroup $X$ of $G$ contains a cyclic subgroup of index two.
	\end{Lem}
	
	\begin{proof}
		We first show that this is the case when $X$ is a Sylow $2$-subgroup. We know that $H$ is contained in some Sylow $2$-subgroup $L$. Moreover, $(L : H) = 2$ since $(G : H)$ is even, but not divisible by $4$. Finally, since $X$ and $L$ are conjugate, $X$ also contains a cyclic subgroup of index two.
		
		Now, assume that $X$ is any $2$-subgroup of $G$ and let $L$ be a Sylow $2$-subgroup containing $X$. We know that there exists a cyclic subgroup $K \leq L$ of index two. Moreover, $K$ is normal by \Cref{lem:index_two_normal}.
		
		We consider two cases. If $X \leq K$, we have that $X$ itself is cyclic, and in particular it has a cyclic subgroup of index $2$ by \Cref{thm:fundamental_cyclic_groups}
		
		If there exists $t \in X \setminus K$, then $L = \langle K, t\rangle$ since $(L : K) = 2$ is prime. By the second isomorphism theorem, $XK$ is a group. Moreover, $XK$ contains both $t$ and $K$, so $XK = L$. The same theorem shows that
		\[
		2 = (L : K) = (XK : K) = (X : X \cap K),
		\]
		so $X \cap K$ has index two inside $X$, and it is cyclic as a subgroup of (the cyclic group) $K$ by \Cref{thm:fundamental_cyclic_groups}.
	\end{proof}
	
	\section{Basic facts about regular origamis and \texorpdfstring{$\translations(g)$}{t(g)}}
	In this section we prove \Cref{thmA:t(g)}, which is stated more precisely below.
	\begin{Thm}\label{thm:not_regular}
		Let $g \geq 2$. If genus-$g$ regular origamis exist, then there exists an integer $m \geq 1$ such that $m \divides 2(g - 1)$, $3 \ndivides m$, $4 \ndivides m$ and
		\[
		\translations(g) = \frac{2(m + 1)}{m} (g - 1).
		\]
		In this case, every translation surface realizing $\translations(g)$ translations is a regular origami up to the action of $\GL^+(2, \RR)$, and belongs to the stratum $\calH(m^{2(g - 1)/m})$.
		
		Otherwise, if no genus-$g$ regular origamis exist, then $\translations(g) = 2(g - 1)$. In this case, every translation surface attaining $\translations(g)$ translations is a normal cover of a torus with two marked points, each with ramification index $g - 1$, and belongs to the principal stratum $\calH(1^{2g - 2})$.
	\end{Thm}
	\Cref{thm:not_regular} is essentially an extension of an argument by Schlage-Puchta and Weitze-Schmithüsen \cite[Lemma 4]{Puchta_Schmithusen}. 
	
	\begin{Notation}
		When the number $m$ in the previous theorem exists, we will denote $\orderm(g) = m$. Otherwise, we set $\orderm(g) = \infty$.
	\end{Notation}
	
	We first state a few useful lemmas.
	
	\subsection{Riemann--Hurwitz formula}
	
	Schlage-Puchta and Weitze-Schmithüsen used the classical Riemann--Hurwitz formula to show that $\translations(g) \leq 4(g - 1)$ for every $g \geq 2$ \cite[Lemma 4]{Puchta_Schmithusen}. The following is a refinement of their argument:
	
	\begin{Lem}\label{lem:origami_bounds}
		Let $X$ be a translation surface of genus $g \geq 2$ with $s$ singularities. We have the following facts:
		\begin{enumerate}
			\item If $X$ belongs to the $\GL^+(2, \RR)$-orbit of a regular origami, then there exists $m \geq 1$ such that $X \in \calH(m^s)$ and
			\[
			|\trans(X)| = \frac{2(m + 1)}{m}(g - 1).
			\]
			\item If $X$ belongs to the $\GL^+(2, \RR)$-orbit of a nonregular origami, then
			\[
			|\trans(X)| \leq 2(g - 1),
			\]
			with equality if and only if $s = 2g - 2$ and, moreover, $X$ is a regular cover of a torus ramified at exactly two distinct points.
			\item If $X$ does not belong to the $\GL^+(2, \RR)$-orbit of an origami, then
			\[
			|\trans(X)| \leq \frac{4}{3}(g - 1).
			\]
		\end{enumerate}
	\end{Lem}
	
	\begin{proof}
		Let $G = \trans(X)$. Consider the quotient translation surface $Y = X/G$ with genus $h \geq 1$. The covering $p\colon X \to Y$ is regular, and its degree $d$ coincides with $|G|$. Moreover, $p$ is ramified at $k \geq 1$ points $P_1, \dots, P_k \in Y$. Each $P_i$ has a number $s_i \geq 1$ of preimages by $p$. Since the covering is regular, all preimages of $P_i$ share the same ramification index $e_i \geq 2$. We have $d = s_i e_i$ for every $1 \leq i \leq k$. In particular, observe that $X$ is in the $\GL^+(2, \RR)$-orbit of an origami if and only if $Y$ is a torus, that is, if and only if $h = 1$. Moreover, such origami is regular if and only if the covering is also ramified at a single point, that is, $k = 1$.
		
		With this data, the Riemann--Hurwitz formula gives
		\begin{align*}
			2g - 2 &= d(2h - 2)+\sum_{i = 1}^k s_i(e_i - 1)\\
			&=d(2h - 2 + k)-\sum_{i = 1}^k s_i.
		\end{align*}
		
		Since $h \geq 1$ and $k\geq 1$, we have
		$2h-2+k \geq 1$. Thus,
		\[
		d = \cfrac{2g-2 + \sum_{i=1}^k s_i}{2h-2+k}.
		\]
		
		Furthermore, we have
		\[
		2g-2 \geq  s \geq \sum_{i=1}^k s_i,
		\]
		with equality in the right hand side inequality if and only if all the singularities of $Y$ are ramification points.
		
		We distinguish the three cases in the statement.
		\begin{enumerate}
			\item If $h = 1$ and $k = 1$, the translation surface $Y$ has no singularities. Hence, $s = s_1$ and
			\[
			d = 2g - 2 + s.
			\]
			Furthermore, every singularity of $X$ shares the same order
			\[
			m = \frac{2g - 2}{s}
			\]
			so $X$ belongs to the stratum $\calH(m^s)$. With this notation,
			\[
			d = 2g - 2 + \frac{2g - 2}{m} = \frac{2(m + 1)}{m}(g - 1).
			\]
			\item If $h = 1$, but $k \geq 2$, we obtain:
			\[
			d \leq \frac{2g-2+s}{k} \leq 2(g - 1),
			\]
			with equality if and only if $k = 2$ and $s = 2g - 2$.
			\item Finally, if $h \geq 2$, we have
			\[
			d \leq \frac{2g-2+s}{2 + k} \leq \frac{2g-2+s}{3} \leq \frac{4}{3}(g - 1).
			\]
		\end{enumerate}
	\end{proof}
	
	We now show that, for every genus $g\geq 2$, there exists a translation surface achieving the bound in the second case of \Cref{lem:origami_bounds}. 
	
	\begin{Lem}\label{lem:2g-2_translations}
		Let $g \geq 2$. There exists a genus-$g$ translation surface with exactly $2(g - 1)$ translations.
	\end{Lem}
	
	\begin{proof}
		
		We exhibit a genus-$g$ origami $X$ on $d = 4g - 4$ squares with translation group $\ZZ/(2g - 2)\ZZ$. Following Matheus' lecture notes \cites[Definition 5]{Matheus:three_lectures}, define $X$ by the following permutations:
		\begin{align*}
			\sigma_{\mathrm{h}}(i) &= i+1 \bmod d\\
			\sigma_{\mathrm{v}}(i) &=
			\begin{cases}
				2g-2+i \mod{4g - 4} &\text{if $i$ is even} \\
				i &\text{if $i$ is odd}.
			\end{cases}
		\end{align*}
		This one-cylinder origami belongs to the stratum $\mathcal{H}(1^{2g-2})$ and its translation group is exactly $\ZZ/(2g-2)\ZZ$. Indeed, first observe that the covering $X \to \TT$ is not regular, so the number of translations is at most $2(g - 1)$ by \Cref{lem:origami_bounds}. Moreover, for each $k \in \{0, 2, 4, \dotsc, 4g - 6\}$, the map sending the square labeled $i$ to the square labeled $i + k \mod{4g - 4}$ defines a translation. Finally, these translations commute.
		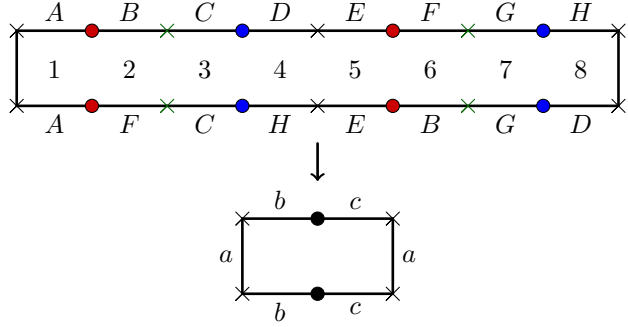
\begin{figure}
			\center
			\definecolor{qqqqff}{rgb}{0,0,1}
			\definecolor{qqwuqq}{rgb}{0,0.39215686274509803,0}
			\definecolor{ccqqqq}{rgb}{0.8,0,0}
			\begin{tikzpicture}[line cap=round,line join=round,x=1cm,y=1cm]
				\clip(-0.5,-3) rectangle (8.5,2);
				\draw [line width=1pt] (0,0)-- (0,1);
				\draw [line width=1pt] (0,1)-- (1,1);
				\draw [line width=1pt] (1,1)-- (2,1);
				\draw [line width=1pt] (2,1)-- (3,1);
				\draw [line width=1pt] (3,1)-- (4,1);
				\draw [line width=1pt] (4,1)-- (5,1);
				\draw [line width=1pt] (5,1)-- (6,1);
				\draw [line width=1pt] (6,1)-- (7,1);
				\draw [line width=1pt] (7,1)-- (8,1);
				\draw [line width=1pt] (8,1)-- (8,0);
				\draw [line width=1pt] (8,0)-- (7,0);
				\draw [line width=1pt] (7,0)-- (6,0);
				\draw [line width=1pt] (6,0)-- (5,0);
				\draw [line width=1pt] (5,0)-- (4,0);
				\draw [line width=1pt] (4,0)-- (3,0);
				\draw [line width=1pt] (3,0)-- (2,0);
				\draw [line width=1pt] (2,0)-- (1,0);
				\draw [line width=1pt] (1,0)-- (0,0);
				\draw [line width=1pt,-to] (4,-0.5) -- (4,-1);
				\draw [line width=1pt] (3,-1.5)-- (4,-1.5);
				\draw [line width=1pt] (4,-1.5)-- (5,-1.5);
				\draw [line width=1pt] (5,-1.5)-- (5,-2.5);
				\draw [line width=1pt] (5,-2.5)-- (4,-2.5);
				\draw [line width=1pt] (4,-2.5)-- (3,-2.5);
				\draw [line width=1pt] (3,-2.5)-- (3,-1.5);
				\draw (0.5,1) node[above] {$A$};
				\draw (0.5,0) node[below] {$A$};
				\draw (1.5,1) node[above] {$B$};
				\draw (5.5,0) node[below] {$B$};
				\draw (2.5,1) node[above] {$C$};
				\draw (2.5,0) node[below] {$C$};
				\draw (3.5,1) node[above] {$D$};
				\draw (7.5,0) node[below] {$D$};
				\draw (4.5,1) node[above] {$E$};
				\draw (4.5,0) node[below] {$E$};
				\draw (5.5,1) node[above] {$F$};
				\draw (1.5,0) node[below] {$F$};
				\draw (6.5,1) node[above] {$G$};
				\draw (6.5,0) node[below] {$G$};
				\draw (7.5,1) node[above] {$H$};
				\draw (3.5,0) node[below] {$H$};
				\draw (3,-2) node[left] {$a$};
				\draw (5,-2) node[right] {$a$};
				\draw (3.5,-1.5) node[above] {$b$};
				\draw (3.5,-2.5) node[below] {$b$};
				\draw (4.5,-1.5) node[above] {$c$};
				\draw (4.5,-2.5) node[below] {$c$};
				\draw (0.5,0.25) node[above] {$1$};
				\draw (1.5,0.25) node[above] {$2$};
				\draw (2.5,0.25) node[above] {$3$};
				\draw (3.5,0.25) node[above] {$4$};
				\draw (4.5,0.25) node[above] {$5$};
				\draw (5.5,0.25) node[above] {$6$};
				\draw (6.5,0.25) node[above] {$7$};
				\draw (7.5,0.25) node[above] {$8$};
				\begin{scriptsize}
					\draw [color=black] (0,0)-- ++(-2.5pt,-2.5pt) -- ++(5pt,5pt) ++(-5pt,0) -- ++(5pt,-5pt);
					\draw [color=black] (0,1)-- ++(-2.5pt,-2.5pt) -- ++(5pt,5pt) ++(-5pt,0) -- ++(5pt,-5pt);
					\draw [fill=ccqqqq] (1,1) circle (2.5pt);
					\draw [color=qqwuqq] (2,1)-- ++(-2.5pt,-2.5pt) -- ++(5pt,5pt) ++(-5pt,0) -- ++(5pt,-5pt);
					\draw [fill=qqqqff] (3,1) circle (2.5pt);
					\draw [color=black] (4,1)-- ++(-2.5pt,-2.5pt) -- ++(5pt,5pt) ++(-5pt,0) -- ++(5pt,-5pt);
					\draw [fill=ccqqqq] (5,1) circle (2.5pt);
					\draw [color=qqwuqq] (6,1)-- ++(-2.5pt,-2.5pt) -- ++(5pt,5pt) ++(-5pt,0) -- ++(5pt,-5pt);
					\draw [fill=qqqqff] (7,1) circle (2.5pt);
					\draw [color=black] (8,1)-- ++(-2.5pt,-2.5pt) -- ++(5pt,5pt) ++(-5pt,0) -- ++(5pt,-5pt);
					\draw [color=black] (8,0)-- ++(-2.5pt,-2.5pt) -- ++(5pt,5pt) ++(-5pt,0) -- ++(5pt,-5pt);
					\draw [fill=qqqqff] (7,0) circle (2.5pt);
					\draw [color=qqwuqq] (6,0)-- ++(-2.5pt,-2.5pt) -- ++(5pt,5pt) ++(-5pt,0) -- ++(5pt,-5pt);
					\draw [fill=ccqqqq] (5,0) circle (2.5pt);
					\draw [color=black] (4,0)-- ++(-2.5pt,-2.5pt) -- ++(5pt,5pt) ++(-5pt,0) -- ++(5pt,-5pt);
					\draw [fill=qqqqff] (3,0) circle (2.5pt);
					\draw [color=qqwuqq] (2,0)-- ++(-2.5pt,-2.5pt) -- ++(5pt,5pt) ++(-5pt,0) -- ++(5pt,-5pt);
					\draw [fill=ccqqqq] (1,0) circle (2.5pt);
					\draw [color=black] (3,-1.5)-- ++(-2.5pt,-2.5pt) -- ++(5pt,5pt) ++(-5pt,0) -- ++(5pt,-5pt);
					\draw [fill=black] (4,-1.5) circle (2.5pt);
					\draw [color=black] (5,-1.5)-- ++(-2.5pt,-2.5pt) -- ++(5pt,5pt) ++(-5pt,0) -- ++(5pt,-5pt);
					\draw [color=black] (5,-2.5)-- ++(-2.5pt,-2.5pt) -- ++(5pt,5pt) ++(-5pt,0) -- ++(5pt,-5pt);
					\draw [fill=black] (4,-2.5) circle (2.5pt);
					\draw [color=black] (3,-2.5)-- ++(-2.5pt,-2.5pt) -- ++(5pt,5pt) ++(-5pt,0) -- ++(5pt,-5pt);
				\end{scriptsize}
			\end{tikzpicture}
			\caption{An origami in the stratum $\mathcal{H}_3(1^4)$ constructed as a regular cover over a torus with two marked points and possessing translation group $\ZZ/4\ZZ$. }
			\label{fig:cover_Z/4Z}
		\end{figure}
	\end{proof}
	
	\subsection{Proof of Theorem~\ref{thm:not_regular}}
	We now have all the ingredients for the proof.
	
	\begin{proof}[Proof of \Cref{thm:not_regular}]
		Assume first that genus-$g$ regular origamis exist. By examining the three cases in \Cref{lem:origami_bounds}, we deduce that regular origamis possess strictly more translations than nonregular origamis and nonorigamis. Hence, there exists $m \geq 1$ such that
		\begin{equation}\label{eq:optimal_m}
			\translations(g) = \frac{2(m + 1)}{m}(g - 1),
		\end{equation}
		where $m$ is the order of the singularities of a translation surface $X$ attaining $\translations(g)$ translations. This surface must lie in the $\GL^+(2, \RR)$-orbit of a regular origami in the stratum $\calH(m^s)$, where $s = (2g - 2) / m$. In particular, $m \divides 2(g - 1)$.
		
		If $3 \divides m$, we deduce that $3 \divides (g - 1)$, so there exists a genus-$g$ origami with translation group of order $4(g - 1)$ \cite[Theorem 1]{Puchta_Schmithusen}. This contradicts \Cref{eq:optimal_m}. Similarly, if $4 \divides m$, we deduce that $2 \divides (g - 1)$ and arrive at a similar contradiction.
		
		If no genus-$g$ origami is regular, we combine the last two cases in \Cref{lem:origami_bounds} to obtain that $\translations(g) \leq 2(g - 1)$. By \Cref{lem:2g-2_translations}, there exists a genus-$g$ origami with this number of translations, and such a surface must lie in the $\GL^+(2, \RR)$-orbit of an origami covering a torus and ramified at exactly two distinct points.
	\end{proof}
	
	\subsection{Building regular origamis}\label{sec:buiding_regular_origamis} Given a finite group $G$ and elements $x, y \in G$ such that $G = \langle x, y\rangle$, we can build a regular origami whose squares are labeled by the elements of $G$, and whose horizontal and vertical permutations are given by (left) multiplication by $G$. If $H$ is the cyclic group generated by $[x, y]$, the resulting regular origami belongs to the stratum $\calH(m^s)$, where $m = |H| - 1$ and $s = (G : H)$. Since $sm = 2g - 2$, we get
	\[
	|G| = |H| \cdot (G : H) = (m + 1) s = 2g - 2 + \frac{2g - 2}{m} = \frac{2(m + 1)}{m} (g - 1).
	\]
	
	Equivalently, a regular origami exists in the stratum $\calH(m^s)$ if and only if there exists a group $G$ or order $(2(m + 1)/m)(g - 1)$, together with two generators $x, y \in G$ such that $H = \langle[x, y]\rangle$ has order $m + 1$ and index $s$ \cite[Remark 2.9]{Flake_Thevis}.
	
	\subsubsection{Direct products} A useful tool to build a regular origami inside a stratum of a prescribed form is to use a direct product between a known translation group $G$ and a cyclic group $\ZZ/k\ZZ$. Concretely, we generalize some ideas of Schlage-Puchta and Weitze-Schmithüsen \cite[Proposition 11]{Puchta_Schmithusen} to include the case where $|G|$ and $k$ are possibly not coprime.
	
	\begin{Lem}\label{lem:extension_cyclic}
		Let $G$ be a group generated by two elements $x, y \in G$ of orders $\alpha$ and $\beta$. Let $k$ be coprime with $\gcd(\alpha, \beta)$. Consider the group $H = G \times \ZZ / k\ZZ$. Then, there exist elements $a, b \in H$ such that $H = \langle a, b\rangle$ and $\ord([a, b]) = \ord([x, y])$.
		
		Furthermore, when $k$ is coprime to $\alpha$, one can choose $a = (x,1)$ and $b=(y,0)$. 
		
		Finally, if the regular origami induced by the group $G$ and the generators $x, y \in G$ lies in the stratum $\calH(m^s)$, then the new regular origami induced by the group $H$ and the generators $a, b \in H$ lies in the stratum $\calH(m^{ks})$.
	\end{Lem}
	
	\begin{proof}
		Since $k$ is coprime with $\gcd(\alpha, \beta)$, we can write $k = ts$ where:
		\begin{itemize}
			\item $t$ and $s$ are coprime;
			\item $t$ is coprime with $\alpha$; and
			\item $s$ is coprime with $\beta$.
		\end{itemize}
		Now, since $t$ and $s$ are coprime, there is an isomorphism
		\[
		H = G \times \ZZ/k\ZZ \simeq G \times \ZZ/t\ZZ \times \ZZ/s\ZZ.
		\]
		Take the elements $a = (x, 1,0)$ and $b = (y, 0,1)$ of $H$.
		Since $\ZZ/k\ZZ$ is abelian, we have $[a, b] = ([x, y], 0, 0)$ and, therefore, $\ord([a, b]) = \ord([x, y])$.
		
		Furthermore, since $t$ is coprime with $\alpha$, there exists $u$ such that $u \alpha \equiv -1 \mod{t}$, and we have that
		\[
		a^{u \alpha + 1} = (x^{u \alpha + 1}, u \alpha + 1, 0) = (x, 0, 0),
		\]
		so $(x,0,0) \in \langle a,b \rangle$. Moreover, if $e \in G$ denotes the identity element of $G$, we have
		\[
		(e, 1, 0) = (x^{\alpha-1}, 0, 0)(x, 1, 0) \in \langle a,b \rangle.
		\]
		Analogously, we obtain that
		$(y,0,0) \in \langle a,b \rangle$ and $(e,0,1) \in \langle a,b \rangle$. Therefore, we deduce that $H \supseteq G \times \{1\} \times \{1\}$, $H \supseteq \{e\} \times \ZZ/t\ZZ \times \{1\}$, and $H \supseteq \{e\} \times \{1\} \times \ZZ/s\ZZ$, so $H =\langle a,b \rangle$. 
		
		Finally, the orders of $[x, y]$ in $G$ and $[a, b]$ in $H$ match, and the index of $\langle[a, b]\rangle$ inside $H$ is $ks$, so the discussion at the beginning of \Cref{sec:buiding_regular_origamis} shows that the regular origami induced by $H$ and the generators $a, b \in H$ lies in the stratum $\calH(m^{ks})$.
	\end{proof}
	
	\subsubsection{Examples of regular origamis.}\label{sec:examples_origamis}
	We know provide a few examples of regular origamis that lie in specific strata, for later use.
	
	\begin{enumerate}[wide,itemsep=1ex]
		\item For any integer $k \geq 1$, the stratum $\calH(k^4)$ contains regular origamis. Indeed, consider the group $\D_{4(k+1)} \simeq \ZZ/2(k+1)\ZZ \rtimes_{-1} \ZZ/2\ZZ$, together with the generators $x = (1,0)$ and $y= (0,1)$. We have $[x,y] = [y,x]^{-1} = x^2$; it has order $k+1$ and generates a group of index $4$. Thus, the induced regular origami lies in $\calH(k^4)$. \label{i:example_k^4}
		\item For any even integers $k, \ell \geq 2$, the stratum $\calH(k^\ell)$ contains regular origamis. Indeed, first consider the group $\D_{2(k+1)} = \ZZ/(k+1)\ZZ \rtimes_{-1}\ZZ/2\ZZ$, together with the generators $x=(1,0)$ and $y=(0,1)$. We have $[x, y] = [y, x]^{-1} = x^2$; it has order $k + 1$ (since $k$ is even) and generates a group of index $2$. Thus, the induced regular origami lies in $\calH(k^2)$. 
		
		Now, observe that $\gcd(\ord(x), \ord(y)) = \gcd(k + 1, 2) = 1$. Take $\lambda = \ell/2$ and apply \Cref{lem:extension_cyclic} using the groups $\D_{2(k+1)}$ and $\ZZ/\lambda\ZZ$, together with the generators $x, y \in \D_{2(k+1)}$. The resulting regular origami lies in the stratum $\calH(k^\ell)$.
		\label{i:example_2k^2l}
		
		\item For any $g \geq 2$ with $9 \divides (g - 1)$, the stratum $\calH(2^{g-1})$ contains regular origamis. Indeed, first consider the group $\M_{3^{\alpha+1}} = \ZZ/3^\alpha \ZZ \rtimes_{3^{\alpha-1}+1} \ZZ/3\ZZ$, together with the generators $x=(1,0)$ and $y = (0,1)$. We have $[x, y] = [y, x]^{-1} = x^{3^{\alpha - 1}}$; it has order $3$ and generates a group of index $9$. Thus, the induced regular origami lies in the stratum $\calH(2^s)$, for $s = 3^\alpha$.
		
		Now, observe that $\gcd(\ord(x), \ord(y)) = \gcd(3^\alpha, 3) = 3$. Write $g - 1 = 3^{\alpha} \lambda$, for $\alpha \geq 2$ and $\lambda$ not divisible by $3$. We use \Cref{lem:extension_cyclic} with the groups $\M_{3^\alpha}$ and $\ZZ/\lambda\ZZ$, together with the generators $x, y \in \M_{3^\alpha}$. The resulting regular origami lies in the stratum $\calH(2^{g-1})$.
		\label{i:example_M_3^alpha}
	\end{enumerate}
	
	In fact, it is possible to completely characterize the strata where regular origamis with a translation group isomorphic to a semidirect product of two cyclic groups. See \Cref{prop:derived_sgp_semi_direct_product}.
	
	\section{The sets \texorpdfstring{$\genera(2)$}{G(2)} and \texorpdfstring{$\genera(2^\alpha - 1)$}{G(2\textasciicircum alpha - 1)} for \texorpdfstring{$\alpha \geq 1$}{alpha >= 2} are empty}\label{sec:G(2)_and_G(2^a)}
	
	Now that we have proven \cref{thmA:t(g)}, we study the set 
	\[ \genera(m) = \left\{ g \geq 2 \st \translations(g) = \frac{2(m+1)}{m}(g-1)\right\}\]
	In this section, we study some particular values of $m$, namely $m = 2$ and $m$ of the form $2^\alpha - 1$ for $\alpha \geq 1$. Using elementary group-theoretic methods, we show that $\genera(m)$ is empty for such values of $m$, proving \Cref{thm:G_m_empty}.
	
	\subsection{The set \texorpdfstring{$\genera(2)$}{G(2)}} We start with the case of $m = 2$. This result will be crucial in later sections to compute certain values of $\orderm(g)$.
	
	\begin{Thm}\label{thm:G(2)_empty}
		The set $\genera(2)$ is empty, that is, $\translations(g) \neq 3(g - 1)$ for every $g \geq 2$.
	\end{Thm}
	
	In fact, we completely classify the set of genera such that the stratum $\calH(2^{g-1})$ contains regular origamis.
	\begin{Thm}\label{thm:H(2^g-1)}
		There exist regular origamis in $\calH(2^{g-1})$ if and only if $g-1$ is even or $9 \divides (g-1)$.
	\end{Thm}
	Since all the genera in the statement of \Cref{thm:H(2^g-1)} belong to $\genera(1)$, this shows that regular origamis in $\calH(2^{g-1})$ never achieve $\translations(g)$, proving \Cref{thm:G(2)_empty}.
	
	\begin{proof}[Proof of \Cref{thm:H(2^g-1)}]
		From the discussion of \Cref{sec:buiding_regular_origamis}, the statement is equivalent to the existence of a group $G$ of order $n = 3(g-1)$ generated by two elements whose commutator has order $3$. If $g-1$ is even or if $9 \divides (g-1)$, we already know that a regular origami exists in $\calH(2^{g-1})$. More precisely, if $g-1$ is even, this is covered in Example~\eqref{i:example_2k^2l}; if $9 \divides (g - 1)$, this is covered in Example~\eqref{i:example_M_3^alpha}.
		We will therefore assume that $g-1$ is odd and that $9 \ndivides (g-1)$.
		
		Observe that $2 \ndivides n$. Thus, by the Feit--Thompson theorem (\Cref{thm:feit-thompson}), $G$ is solvable. As a consequence, since $3 \divides n$, we know from \Cref{thm:Hall} that there exists a $3'$-Hall subgroup $U$ of $G$.
		
		From $9 \ndivides (g-1)$, we obtain that $27 \ndivides n$, and therefore that $(G : U) \divides 9$. Let $\ell = n/(G:U) = |U|$. By construction, $3 \ndivides \ell$. We will show that $U$ is normal.
		
		First, the number of conjugates of $U$ is $k = (G : \N_G(U))$. But,
		\[
		(G : U) = (G : \N_G(U))(\N_G(U) : U),
		\]
		and, therefore, $k=(G : \N_G(U))$ must be either $1, 3$ or $9$. 
		
		Using the action of $U$ on the set $\Omega = \{gUg^{-1} \st g \in G\}$, Schlage-Puchta and Weitze-Schmithüsen show \cite[Lemma 15]{Puchta_Schmithusen} that there exist integers $a_p \geq 0$ for each prime divisor $p$ of $\ell$ such that: 
		\[
		k = 1 + \sum_{p \divides \ell} a_p p.
		\]
		Now, since $\ell$ is not divisible by $p = 2$ or $p = 3$, these factors do not appear in the sum. Hence, $k$ cannot be $3$ or $9$. Thus, $k=1$, and $U$ is normal.
		
		As a consequence, the quotient $G/U$ is a group of order either $3$ or $9$, hence it is abelian. Thus, the commutator subgroup of $G$ is a subgroup of $U$, and any commutator has an order dividing $\ell=|U|$. As $3 \ndivides \ell$, the order of any commutator is not $3$.
	\end{proof}
	
	\subsection{The set \texorpdfstring{$\genera(2^\alpha - 1)$}{G(2\textasciicircum alpha - 1)}}
	We now turn our attention to $m = 2^\alpha - 1$.
	
	We will prove the following:
	
	\begin{Thm}\label{thm:G(2^k-1)_empty}
		The set $\genera(2^\alpha - 1)$ is empty for each $\alpha \geq 2$.
	\end{Thm}
	
	We again work in the group theoretic setting. We will show a stronger version of this:
	
	\begin{Thm}\label{thm:G(2^k-1)_empty_group}
		Let $g \geq 2$ be an integer with $g \notin \genera(1)$. Let $\alpha \geq 2$ be an integer. Then, a group $G$ of order
		\[
		n = \frac{2^{\alpha + 1}}{2^\alpha - 1}(g - 1),
		\]
		together with two generators $x, y \in G$ such that $[x, y]$ has order $2^\alpha$, does not exist.
	\end{Thm}
	
	\begin{Rema}
		In geometric terms, this means the stratum $\calH((2^\alpha - 1)^\ell)$ contains no regular origamis if $g = \ell(2^\alpha - 1)/2 + 1 \notin \genera(1)$, that is, if $\ell \equiv 2 \mod{4}$ and $3 \ndivides \ell(2^\alpha - 1)$. Nevertheless, such origamis do exist if $g \in \genera(1)$ \cite[Theorem A]{Flake_Thevis}.
	\end{Rema}
	
	The crux of the proof is the following lemma:
	
	\begin{Lem}\label{lem:lem:G(2^k-1)_normal_2'_Hall}
		In the context of \Cref{thm:G(2^k-1)_empty_group}, $G$ contains a normal $2'$-Hall subgroup.
	\end{Lem}
	
	\begin{proof}
		This is an application of Frobenius $p$-complement theorem (\Cref{thm:Frobenius_p_complement}) together with \Cref{lem:cyclic_index_two}. Indeed, given a $2$-subgroup $X \leq G$, we will show that $\N_G(X)/\C_G(X)$ is a $2$-group.
		
		Let $q > 3$ be prime. Since $H = \langle [x, y]\rangle$ is a cyclic $2$-subgroup of $G$ of index
		\[
		(G : H) = \frac{2}{2^\alpha - 1}(g - 1)
		\]
		and $g - 1$ is odd as $g \notin \genera(1)$, \Cref{lem:cyclic_index_two} shows that $X$ contains an index-two cyclic subgroup and, in particular, that it is isomorphic to one of the groups in \Cref{thm:Burnside}. By \Cref{lem:Burnside_Aut}, $q$ does not divide $|\Aut(X)|$.
		In particular, from \Cref{rk:normalizer_centralizer_automorphism} $q$ does not divide $|\N_G(X) / \C_G(X)|$ either.
		
		On the other hand, we have that $3$ does not divide $|G|$, so it does not divide $|\N_G(X)|$. Thus, the only prime factor of $|\N_G(X) / \C_G(X)|$ is $2$, so this group is a $2$-group. We conclude using the Frobenius $p$-complement theorem.
	\end{proof}
	
	We can now finish the proof of \Cref{thm:G(2^k-1)_empty_group}:
	
	\begin{proof}[Proof of \Cref{thm:G(2^k-1)_empty_group}]
		Let $N$ be a normal $2'$-Hall subgroup of $G$, which exists by \Cref{lem:lem:G(2^k-1)_normal_2'_Hall}.
		
		Let $H = \langle [x, y]\rangle$. Since $2^2 \ndivides (G : H)$, we have that $2^2 \ndivides (G : G')$. Thus, \Cref{lem:p_Sylow_order_p_Hall} shows that
		\[
		G \simeq N \rtimes \ZZ/2\ZZ,
		\]
		with $|N|$ odd. Hence, $2^{\alpha}$ does not divide $|G|$, so no commutator has order $2^{\alpha}$.
	\end{proof}
	
	\section{Regular origamis with translation group \texorpdfstring{$\PSL(2, p)$}{PSL(2, p)}} \label{sec:PSL}
	
	In this section, we prove a more precise version of \Cref{thm:G_m_intro}, namely:
	\begin{Thm}\label{thm:G_m_general}
		Let $m \geq 5$ be prime and assume that $3 \divides (m + 1)$. There exist infinitely many prime numbers $p$ such that, for every $k \geq 1$ not divisible by any prime number $q < m$, we have
		\[
		\frac{k m p(p - 1)(p + 1)}{4(m + 1)} + 1 \in \genera(m).
		\]
		In particular, $\genera(m)$ contains infinitely long arithmetic progressions.
	\end{Thm}
	For this, we will construct regular origamis with translation group of the form $\PSL(2,p) \times \ZZ/k\ZZ$ for suitable prime $p$, integer $k$, and appropriate pairs of generators $(A,B)$, and we will show that they have the largest automorphism group for their genus. The main technical device for this section is producing generators $A$, $B$ of $\SL(2, p)$ such that $[A, B]$ has a desired order $d$ and such that $B$ has order two:
	\begin{Prop}
		\label{prop:generates}
		Let $p > 13$ be a prime number. Then, for any integer $d\geq 6$ satisfying $p\equiv \pm 1 \mod d$, there exists $A \in \mathrm{SL}(2, p)$ such that $A$ and $B = \Big(\begin{smallmatrix} 0 & -1 \\ 1 & 0 \end{smallmatrix} \Big)$ generate $\mathrm{SL}(2,p)$, and $[A,B]$ has order $d$ in $\mathrm{SL}(2, p)$.
	\end{Prop}
	\Cref{prop:generates} is very similar to the work of McCullough and Wanderley~\cite[Theorem 2.2]{McCullough_Wanderley_commutators}: they show that any element can be realized as the commutator of a generating pair $(A,B)$, except for $-\Id$ and those of trace $2$. In fact, their result is sufficient to prove a slightly weaker version of \Cref{thm:G_m_general}, where we additionally assume that $k$ is coprime with the order of $\PSL(2,p)$ (that is, $k$ is coprime with $p$, $p-1$, and $p+1$). This additional assumption makes $k$ depend not only on $m$, but also on $p$.
	To rule out the dependence on $p$, we will additionally need to control $\gcd(\ord(A), \ord(B))$ in order to apply \Cref{lem:extension_cyclic}, and therefore we include a complete (and different) proof.
	
	Before proving this \Cref{prop:generates}, we need a result about generating pairs of $\SL(2, p)$ and two lemmas. Our first lemma shows that order-$d$ elements exist in $\SL(2, p)$ for suitable $d$. This result is well-known and it is actually more general \cite[Theorem 2.4]{Glover_Sjerve_PSL}. Nevertheless, we include a short proof for the sake of completeness.
	
	\begin{Lem}
		\label{lem:SL_order}
		Let $p$ be an odd prime number. Let $d \geq 3$ be such that $p \equiv \pm 1 \mod{d}$. Then, there exists an order-$d$ element $M \in \SL(2, p)$. Moreover, $\tr(M) \neq \pm 2$.
	\end{Lem}
	
	\begin{proof}
		We will construct $M$ explicitly. We consider two cases:
		\begin{itemize}
			\item If $p\equiv 1 \mod d$, we choose $\lambda\in \mathbb{F}^\times_p$ of multiplicative order $d$. This is possible since $\FF_p^\times$ is cyclic of order $p - 1$. The element $M=\diag(\lambda,\lambda^{-1}) \in \mathrm{SL}(2,p)$ has order $d$. 
			
			\item If $p\equiv -1 \mod d$, we choose $\lambda\in \mathbb{F}^\times_{p^2} \setminus \FF_p^\times$ of multiplicative order $d$. This is possible since $\FF_{p^2}^\times$ is cyclic of order $(p - 1)(p + 1)$. Since $\lambda$ is sent to $\lambda^{-1}$ by the unique nontrivial element of the Galois group $\Gal(\FF_{p^2} / \FF_p)$, the element $t=\lambda+\lambda^{-1}$ belongs to $\mathbb{F}_p$. Moreover, $\lambda \neq \lambda^{-1}$ since, otherwise, $\lambda \in \FF_p$.
			
			Now, $\lambda$ and $\lambda^{-1}$ are the roots of the quadratic polynomial $x^2-tx+1$. Take $M=\left( \begin{smallmatrix} t & -1 \\ 1 & 0 \end{smallmatrix} \right) \in \mathrm{SL}(2,p)$. This matrix has characteristic polynomial $x^2-tx+1$, so it is conjugate to $\diag(\lambda, \lambda^{-1})$. Thus, it has order $d$, as required. \qedhere
		\end{itemize}
		In both cases, $\tr(M) \neq \pm 2$ since, otherwise, $\lambda = \pm 1$, which does not have order $d$.
	\end{proof}
	
	We continue with a folklore lemma in elementary number theory:
	\begin{Lem} \label{lem:sum_of_two_squares}
		Let $p \geq 17$ be prime and let $a \in \FF_p^\times$. Then, $a$ can be written as the sum of two nonzero squares in at least two different ways. That is, there exist $s_1, t_1, s_2, t_2 \in \FF_p^\times$ such that $a = s_1^2 + t_1^2 = s_2^2 + t_2^2$, and $\{s_1^2, t_1^2\} \cap \{s_2^2, t_2^2\} = \varnothing$.
	\end{Lem}
	
	\begin{proof}
		Fix $a \in \FF_p^\times$. Consider $Q \colon \FF_p^2 \to \FF_p$ given by $Q(x, y) = x^2 + y^2$ and define $R = \{x^2 \st x \in \FF_p\}$. We have that $|R| = (p + 1)/2$, so $R \cap (a - R) \neq \varnothing$. We get $x^2 = a - y^2$ for some $x, y \in \FF_p$, so $Q(x, y) = a$.
		
		Now, the group $\SO(Q)$ has order at least $p - 1 \geq 16$ \cite[p.\ 141]{Taylor_classical_groups}. Define $S = \SO(Q) \cdot (x, y)$. Since the stabilizer of $(x, y)$ is trivial, we obtain that $|S| \geq 16$.
		
		Finally, $S$ contains at most two elements of the form $(\pm r, 0)$ and at most two elements of the form $(0, \pm r)$, so at least twelve of its elements belong to $(\FF_p^\times)^2$. If $(s_1, t_1)$ is such an element, $S$ contains at most four elements of the form $(\pm s_1, \pm t_1)$ and at most four elements of the form $(\pm t_1, \pm s_1)$, so there exist $(s_2, t_2) \in S$ not of these forms. We obtain that $a = s_1^2 + t_1^2 = s_2^2 + t_2^2$ and $\{s_1^2, t_1^2\} \cap \{s_2^2, t_2^2\} = \varnothing$.
	\end{proof}
	
	Our final ingredient to prove \Cref{thm:G_m_general} is the following number-theoretic fact:
	
	\begin{Lem}\label{lem:infinitely_many_primes}
		Let $m$ be prime with $3 \divides (m + 1)$. Then, there exist infinitely many prime numbers $p$ with $p \equiv \pm 1 \mod{2(m + 1)}$ such that
		\[
		z = \frac{(p - 1)(p + 1)}{4(m + 1)}
		\]
		is an integer and is not divisible by any prime number $q < m$.
	\end{Lem}
	
	\begin{proof}
		We will use the Chinese remainder theorem to show that this holds if $p$ belongs to certain residue classes. Then, infinitely many such prime numbers will exist due to Dirichlet's theorem.
		
		Observe that $z$ is an integer if $p \equiv \pm 1 \mod{2(m + 1)}$. Indeed, this means that $2(m + 1)$ divides either $p - 1$ or $p + 1$, and the other factor of $(p - 1)(p + 1)$ is even.
		
		Let $2^\alpha$, with $\alpha \geq 1$, be the largest power of $2$ that divides $m + 1$. Observe that the largest power of $2$ that divides $4(m + 1)$ is $2^{\alpha + 2}$. For $z$ to be an odd integer, we need the largest power of $2$ that divides $(p - 1)(p + 1)$ also to be $2^{\alpha + 2}$. This is equivalent to
		\[
		p \equiv a \mod{2^{\alpha + 2}},
		\]
		where $a \not\equiv \pm1 \mod{2^{\alpha + 2}}$,
		and $a \equiv \pm 1 \mod{2^{\alpha + 1}}$. Indeed, this equation imposes that the largest power of $2$ that divides $p - 1$ or $p + 1$ is $2^{\alpha + 1}$, and the other factor in $(p - 1)(p + 1)$ is always even (and cannot be divisible by $4$). Consequently, we take $a = 2^{\alpha + 1} \pm 1$.
		
		Now, let $\Pi$ be the set of primes $q$ with $3 \leq q < m$. If $q \in \Pi$, consider the largest power $q^{\beta_q}$ that divides $m + 1$ (possibly, $\beta_q = 0$). For $z$ to be an integer, we need $p \equiv \pm 1 \mod{q^{\beta_q}}$.
		
		Furthermore, observe that $q \ndivides z$ is equivalent to $p \not\equiv \pm 1 \mod{q^{\beta_q + 1}}$. This is equivalent to
		\[
		p \equiv b_q \mod{q^{\beta_q + 1}},
		\]
		where $b_q \not\equiv \pm 1 \mod{q^{\beta_q + 1}}$, and $b_q \equiv \pm 1 \mod{q^{\beta_q}}$. Hence, we take $b_q = \ell_q q^{\beta_q} \pm 1$, where $0 < \ell_q < q$ if $\beta_q > 0$, and $1 < \ell_q < q - 1$ otherwise. The previous argument works for $q = 3$ since $\beta_q \geq 1$ as $3 \divides (m + 1)$ by hypothesis (this is necessary since every prime $p \geq 5$ satisfies $(p - 1)(p + 1) \equiv 0 \mod{3}$, so we need $\beta_q + 1 > 1$).
		
		Combining this information, $p$ must satisfy the system of congruences:
		\begin{align*}
			p &\equiv \pm 1 \mod{2(m + 1)} \\
			p &\equiv 2^{\alpha + 1} \pm 1 \mod{2^{\alpha + 2}} \\
			p &\equiv \ell_q q^{\beta_q} \pm 1 \mod{q^{\beta_q + 1}} \text{ for every $q \in \Pi$},
		\end{align*}
		
		Then, the (generalized) Chinese remainder theorem shows that the system admits a solution if and only if the equations are pairwise compatible modulo the corresponding greatest common divisors. We have:
		\begin{align*}
			\gcd(2(m + 1), 2^{\alpha + 2}) &= 2^{\alpha + 1} \\
			\gcd(2(m + 1), q^{\beta_q + 1}) &= q^{\beta_q} \\
			\gcd(2^{\alpha + 2}, q^{\beta_q + 1}) &= 1.
		\end{align*}
		Thus, the compatibility is automatic as long as the choice of signs is consistent (that is, all signs are ``$+1$'', or all signs are ``$-1$''), independently of the $\ell_q$.
		
		Finally, observe that
		\[
		\lcm(\{2(m + 1), 2^{\alpha + 2}\} \cup \{q^{\beta_q + 1} \st q \in \Pi\}) = 4(m + 1)Q,
		\]
		where $Q$ is the product of all $q \in \Pi$. Hence, the Chinese remainder theorem also shows that there exists $t$ such that $p$ is a solution if $p \equiv t \mod{4(m + 1)Q}$. Since $p$ is coprime with $4(m + 1)Q$ by construction, so is $t$.
		
		We conclude by Dirichlet's theorem: infinitely many primes belong to this residue class $t$ modulo $4(m + 1)Q$.
	\end{proof}
	
	\begin{Rema}
		When $m$ is a Sophie Germain prime, the smallest prime number $p$ as in the previous lemma is $p = 2m + 1$. This choice yields $z = m$.
	\end{Rema}
	
	To show that two matrices $A$ and $B$ generate $\mathrm{SL}(2,p)$, we will use the results of McCullough and Wanderley \cite{McCullough_Wanderley_generators}, building on the results of Macbeath \cite{Macbeath}.
	
	\begin{Thm}[{\cite[Section 11]{McCullough_Wanderley_generators}}]\label{thm:McCullough_Wanderley}
		Let $p \geq 13$ be a prime number. Two elements $A, B \in \SL(2, p)$ are generators if and only if:
		\begin{enumerate}
			\item at least two of the numbers $\tr(A)$, $\tr(B)$ and $\tr(AB)$ are nonzero;
			\label{i:generates_cond_1}
			\item $\tr([A, B]) \neq 2$; and \label{i:generates_cond_2}
			\item $\langle A, B\rangle$ is not isomorphic to $\A_4$, $\Sym_4$ or $\A_5$. \label{i:generates_cond_3}
		\end{enumerate}
	\end{Thm}
	
	The conditions they find are, in fact, more general since they deal with any finite field. Nevertheless, \Cref{thm:McCullough_Wanderley} is enough for our purposes, as we only work with fields of odd-prime order.
	
	Finally, McCullough and Wanderley provide a full classification of conjugacy classes in $\SL(2, p)$ \cite[Proposition 2.3]{McCullough_Wanderley_commutators}. We recall part of it: $A, B \in \SL(2, p)$ with $\tr(A), \tr(B) \neq \pm 2$ are conjugate if and only if $\tr(A) = \tr(B)$.
	
	We can now prove the proposition:

	\begin{proof}[Proof of~\Cref{prop:generates}]
		
		We take $A = \Big(\begin{smallmatrix} a & b \\ c & d \end{smallmatrix} \Big) \in \SL(2, p)$, where $a,b,c,d \in \mathbb{F}_p$ will be chosen appropriately.
		
		Let $M \in \SL(2, p)$ be an element of order $d$, which exists by \Cref{lem:SL_order}. Take $x = \tr(M)$; we have $x \neq \pm 2$. On the one hand, we will ensure that $\tr([A, B]) = x$. This implies $A$ has order $d$ by the classification of conjugacy classes above.
		
		On the other hand, by \Cref{thm:McCullough_Wanderley}, we need to verify conditions \eqref{i:generates_cond_1}, \eqref{i:generates_cond_2}, and \eqref{i:generates_cond_3} to obtain $\SL(2, p) = \langle A, B\rangle$. In fact, it will be enough to focus on condition~\eqref{i:generates_cond_1}. Indeed, condition~\eqref{i:generates_cond_2} follows automatically from $\tr([A, B]) = x$. Moreover, condition~\eqref{i:generates_cond_3} is automatic from the fact that $d \geq 6$, since the only possible orders of commutators in $\A_4$, $\Sym_4$, and $\Sym_5$ are $1$, $2$, $3$, and $5$.
		
		Observe that
		\begin{align*}
			\tr(A) &= a + d \\
			\tr(AB) &= b - c \\
			\tr([A,B]) &= a^2 + b^2 + c^2 + d^2.
		\end{align*}
		Consider the equation
		\begin{equation}\label{eq:SL_t_s}
			x^2 - 4 =(2s)^2 + t^2,
		\end{equation}
		with variables $s$ and $t$. Since $x \neq \pm 2$, we have $x^2 - 4 \neq 0$. By \Cref{lem:sum_of_two_squares}, there exist two pairs of solutions $(s_1, t_1), (s_2, t_2) \in (\FF_p^\times)^2$ of \Cref{eq:SL_t_s} such that $\{s_1^2, t_1^2\} \cap \{s_2^2, t_2^2\} = \varnothing$. If $t_1 \neq \pm x$, we define $(s, t) = (s_1, t_1)$. Otherwise, we have $t_2 \neq \pm x$ since $t_1^2 \neq t_2^2$, and we take $(s, t) = (s_2, t_2)$.
		
		Now, take $u = (x + t) / 2$, which is nonzero by our choice of $t$. A straightforward computation shows that
		\begin{equation}\label{eq:SL_u}
			u x - (u^2 + 1) = s^2.
		\end{equation}
		
		Using \Cref{lem:sum_of_two_squares} again, we continue by taking $(a_1, c_1), (a_2, c_2) \in (\FF_p^\times)^2$ such that $u = a_1^2 + c_1^2 = a_2^2 + c_2^2$ and $\{a_1^2, c_1^2\} \cap \{a_2^2, c_2^2\} = \varnothing$.
		
		Now, for $i \in \{1, 2\}$, we define:
		\[
		b_i^{\pm} = \frac{-c_i \pm a_i s}{a_i^2 + c_i^2}, \qquad d_i^{\pm} = \frac{1 + b_i^{\pm} c_i}{a_i}.
		\]
		
		The choice of $d_i^{\pm}$ directly implies $a_i d_i^{\pm} - b_i^{\pm} c_i = 1$. Moreover,
		\begin{align*}
			a_i^2 + (b_i^{\pm})^2 + c_i^2 + (d_i^{\pm})^2 &= \frac{(a_i^2 + c_i^2)^2 + 1 + s^2}{a_i^2 + c_i^2} \\
			&= \frac{u^2 + 1 + ux - (u^2 + 1)}{u} \\
			& = x,
		\end{align*}
		were we used that $a_i^2 + c_i^2 = u$ and \Cref{eq:SL_u}.
		
		We will show that, for some choice of $i \in \{1, 2\}$ and $\varepsilon \in \{+, -\}$, we have $a_i + d_i^\varepsilon \neq 0$ and $b_i^\varepsilon - c_i \neq 0$, so we can take $a = a_i$, $b = b_i^\varepsilon$, $c = c_i$, and $d = d_i^\varepsilon$ to finish the proof.
		
		Observe that
		\begin{align*}
			a_i + d_i^{\pm} &= \frac{a_i(u + 1) \pm c_i s}{u} \\
			b_i^{\pm} - c_i &= \frac{-c_i(u + 1) \pm a_i s}{u}.
		\end{align*}
		
		If one of these quantities vanishes, we can solve for $s$ to obtain
		\begin{align*}
			a_i + d_i^{\pm} &= 0 \implies s = \mp \frac{a_i}{c_i} (u + 1) \\
			b_i^{\pm} - c_i &= 0 \implies s = \pm \frac{c_i}{a_i} (u + 1).
		\end{align*}
		
		In particular, if any of these numbers vanish, then $u \neq -1$ since $s \neq 0$. Furthermore, we deduce that $a_i + d_i^\varepsilon = 0$ for both choices of $\varepsilon \in \{+, -\}$ implies $s = 0$, which is impossible. Similarly, $b_i^\varepsilon - c_i = 0$ for both choices of $\varepsilon \in \{+, -\}$ is impossible.
		
		Thus, if one of the numbers $a_i + d_i^\varepsilon$ or $b_i^\varepsilon - c_i$ vanishes for both choices of $\varepsilon \in \{+, -\}$, we only have two possible cases: either both $a_i + d_i^+$ and $b_i^- - c_i$ vanish, or both $a_i + d_i^-$ and $b_i^+ - c_i = 0$ do. In both cases, solving for $s$ as above and canceling $u + 1$ out yields $a_i/c_i = c_i/a_i$, so $c_i^2 = a_i^2$.
		
		Finally, we see that $u = a_i^2 + c_i^2$ is written as a sum of two \emph{equal} squares. This can only happen for a single choice of $i \in \{1, 2\}$ since $\{a_1^2, c_1^2\} \cap \{a_2^2, c_2^2\} = \varnothing$. Thus, if $j$ is such that $\{i, j\} = \{1, 2\}$, we deduce that $a_j + d_j^\varepsilon \neq 0$ and $b_j^\varepsilon - c_j \neq 0$ for at least one of the choices of $\varepsilon \in \{+, -\}$, completing the argument.
	\end{proof}
	
	\begin{Rema}
		The only reason why we need $d \geq 6$ in the previous proof is to rule out the groups $\A_4$, $\Sym_4$, and $\A_5$, that is, to establish condition \eqref{i:generates_cond_3} in \Cref{thm:McCullough_Wanderley}. This condition is actually not necessary, since $\langle A, B \rangle$ is one of these groups only in very particular situations. Since we only need the result for $d \geq 6$, we refrain from stating this more general version and refer the reader to the work of McCullough \cite{McCullough_exceptional_groups} for more details.
	\end{Rema}
	
	We now have all the tools to prove the main result of this section.
	\begin{proof}[Proof of ~\Cref{thm:G_m_general}]
		Let $p$ be as in \Cref{lem:infinitely_many_primes}, and let $k$ be as in the statement. We take $G=\PSL(2, p)$. We define the group $H = G \times \ZZ/k\ZZ$, together with two generators $a, b \in H$ with $\ord([a, b]) = m + 1$. Such generators exist.
		Indeed, by \Cref{prop:generates}, since $p \equiv \pm 1 \mod{2(m+1)}$ there exist two generators $x,y\in \SL(2,p)$ whose commutator $[x,y]$ has order $2(m+1)$ in $\SL(2,p)$, or $m+1$ in $G$.
		Moreover, in \Cref{prop:generates}, we can take the order of $y$ to be two. Thus $k$ is coprime with $\gcd(\ord(x),\ord(y))$, and \Cref{lem:extension_cyclic} applies.
		
		Let $n$ be the order of $H$, that is,
		\[
		n = \frac{k p (p - 1)(p + 1)}{2}.
		\]
		Now, consider the regular origami induced by $H$, $a$, and $b$, and let $g$ be its genus. We have that
		\[
		\frac{k p (p - 1)(p + 1)}{2} = n = \frac{2(m + 1)}{m}(g - 1),
		\]
		so
		\[
		g - 1 = \frac{mk p (p - 1)(p + 1)}{4(m + 1)}.
		\]
		As a consequence, observe that our assumptions on $p$ and $k$ guarantee that $g-1$ is not divisible by any prime number $q < m$. In particular $2 \ndivides (g - 1)$ and $3 \ndivides (g - 1)$, so $g \notin \genera(1)$. Thus, $\orderm(g) \in \{ 5, \dots, m \}$ (recall that $\genera(2)$ is empty by \Cref{thm:G(2)_empty}, and $\orderm(g) \neq 3,4$ by \Cref{thm:not_regular}). Further, if $5 \leq d < m$, we see that $d \ndivides 2(g-1)$ and therefore $g \notin \genera(d)$. Hence, $g \in \genera(m)$.
	\end{proof}
	
	\begin{Rema}
		In the particular case of $m = 5$, the previous proof shows that every prime number $p > 13$ satisfying that $p \bmod 72$ is $11$, $13$, $59$, or $61$ allows us to produce an infinitely long arithmetic progression inside $\genera(5)$, where the condition $p > 13$ is only needed because of the use of \Cref{lem:sum_of_two_squares}. Nevertheless, the values $p = 11$ and $p = 13$ also work, although the general proof fails.
		
		To see this, we can exhibit explicit matrices that generate the groups $\PSL(2, 11)$ and $\PSL(2, 13)$, one of which has order $2$, and with a commutator of order $6$.
		
		For example, the following choices work for $\PSL(2, 11)$:
		\[
		A = \begin{pmatrix}
			1 & 2 \\ 0 & 1
		\end{pmatrix} \qquad \text{and} \qquad B = \begin{pmatrix}
			0 & -1 \\ 1 & 0
		\end{pmatrix}.
		\]
		Similarly, these matrices do the trick for $\PSL(2, 13)$:
		\[
		A = \begin{pmatrix}
			2 & 4 \\ 0 & 7
		\end{pmatrix} \qquad \text{and} \qquad B = \begin{pmatrix}
			0 & -1 \\ 1 & 0
		\end{pmatrix}.
		\]
		Primes smaller than $17$ only arise in the case of $m = 5$, so it is not necessary to consider them in other cases.
		
		Finally, in the case of $m = 11$, the condition that $7 \ndivides z$ is not needed in \Cref{lem:infinitely_many_primes}, since we know that $\genera(7)$ is empty (\Cref{thm:G(2^k-1)_empty}). This case is somewhat simplified to requiring only that $2 \ndivides z$, $3 \ndivides z$, and $5 \ndivides z$. We obtain that any prime $p$ such that $p \bmod{720}$ is $23$, $167$, $263$, $313$, $407$, $457$, $553$, $697$ induces an infinitely long arithmetic progression inside $\genera(11)$.
	\end{Rema}
	
	\section{Infinite families of genera in \texorpdfstring{$\genera(\infty)$}{G(inf)}}
	
	In this section, we prove \Cref{thm:no_regular_origamis_intro,thm:SG_intro,thm:all_strata} by producing infinite families of genera $g \geq 2$ with no genus-$g$ regular origamis and studying the special case of Sophie Germain primes (i.e. a prime $p$ such that $2p + 1$ is also prime):
	\begin{Thm}\label{thm:no_regular_origamis}
		Let $g \geq 6$ be of the form:
		\begin{enumerate}
			\item $g = p + 1$, where $p \geq 5$ is prime; \label{i:no_regular_origamis_1}
			\item $g = p^2 + 1$, where $p$ is prime, but it is not a Sophie Germain prime; \label{i:no_regular_origamis_2}
			\item $g = pq + 1$, where $p, q \geq 5$ are distinct primes. \label{i:no_regular_origamis_3}
		\end{enumerate}
		Then there exist no genus-$g$ regular origamis.
		
		On the other hand, if $p \geq 5$ is a Sophie Germain prime there exist regular origamis of genus $p^2+1$. Such origamis belong to the stratum $\calH(2p^p)$ and have translation group $\ZZ/(2p+1)\ZZ \rtimes \ZZ/p\ZZ$. In particular, $p^2+1 \in \genera(2p)$.
	\end{Thm}
	
	\begin{Rema}
		Combining the last part of \Cref{thm:no_regular_origamis} with \Cref{lem:extension_cyclic}, we obtain that for any Sophie Germain prime $p \geq 5$ and any integer $\ell \geq 1$, it is possible to construct a regular origami of genus $g = \ell p^2+1$ with translation group $(\ZZ/(2p+1)\ZZ \rtimes \ZZ/p\ZZ) \times \ZZ/\ell\ZZ$. In particular $\ell p^2 + 1 \notin \genera(\infty)$.
	\end{Rema}
	
	This proof will be done in several steps. The crux of the proof is showing that regular origamis in the strata $\calH(g - 1, g - 1)$, $\calH(2p^q)$, and $\calH(p^{2q})$, for $p$, $q$ prime, can only exist in very particular situations. These cases will be done in \Cref{s:H(g-1 g-1)}, \Cref{s:H(2p^q)}, and \Cref{s:H(p^2q)}, respectively. Recall from \Cref{sec:buiding_regular_origamis} that if $m \divides (2g - 2)$, the existence of a regular origami in $\calH(m^s)$, for $s = (2g - 2) / m$, is equivalent to the existence of a group $G$ of order
	\[
	n = |G| = \frac{2(m + 1)}{m}(g - 1)
	\]
	that is generated by two elements, $x, y \in G$, such that $[x, y]$ has order $m + 1$. To rule out the existence of regular origami in such strata, we will use several group-theoretic tools.
	
	Finally, we will combine these results to complete the proof of \Cref{thm:no_regular_origamis} in \Cref{s:no_regular_origamis_proof}.
	
	\subsection{The stratum \texorpdfstring{$\calH(g - 1, g - 1)$}{H(g - 1, g - 1)}} \label{s:H(g-1 g-1)}
	We will show that $\calH(g - 1, g - 1)$ can only contain regular origamis if $g$ is odd:
	
	\begin{Lem}\label{lem:H(g-1 g-1)}
		Let $g \geq 2$. If the stratum $\mathcal{H}(g-1,g-1)$ contains regular origamis if and only if $g$ is odd.
	\end{Lem}
	
	\begin{proof}
		Assume that regular origamis exist in the stratum $\mathcal{H}(g-1,g-1)$. From the group-theoretic viewpoint, this means assuming the existence of a group $G$ of order $2g$, generated by two elements $x,y \in G$ whose commutator has order $g$. We will show that $g$ is odd.
		
		Consider the cyclic subgroup $H = \langle[x,y]\rangle \leq G'$ of $G$. Since $|H| = g$, we see that $H$ has index two in $G$. Hence, $H$ is normal by \Cref{lem:index_two_normal}.
		
		We now apply \Cref{cor:p_Sylow_order_p_cyclic} with $p = 2$ to deduce that
		\[
		G \simeq \ZZ/g\ZZ \rtimes \ZZ/2\ZZ,
		\]
		with $2 \ndivides g$, so $g$ is odd.
		
		The converse follows directly from the first example of \Cref{sec:examples_origamis}
	\end{proof}
	
	\subsection{The stratum \texorpdfstring{$\mathcal{H}(2k^q)$}{H(2k\textasciicircum q)}} \label{s:H(2p^q)}
	
	If $q$ is an odd prime and $k$ is an integer, regular origamis in the stratum $\calH(2k^q)$ can only exist under very special conditions:
	
	\begin{Thm}\label{thm:stratum_H(2p^q)}
		Let $k \geq 1$ be an integer and let $q$ be an odd prime. There exist regular origamis in the stratum $\mathcal{H}(2k^q)$ if and only if there exists $d \in \{2,3,\dots,2k\}$ such that
		\begin{itemize}
			\item $d^q\equiv 1 \mod 2k+1$;
			\item $d-1$ is coprime with $2k+1$.
		\end{itemize}
		In this case, the translation group of such a regular origami is isomorphic to the semidirect product $\ZZ/(2k+1)\ZZ \rtimes_d \ZZ/q\ZZ$.
		
		Furthermore, the existence of such $d$ is equivalent to all prime factors of $2k+1$ being congruent to $1$ modulo $q$.
	\end{Thm}
	
	\begin{Rema} The case $\calH(k^2)$ reduces to the case $\calH(g-1,g-1)$, which has been dealt with in \Cref{lem:H(g-1 g-1)}. Moreover, since $q$ is odd, $\calH(u^q)$ is empty if $u$ is odd, so we assume $u = 2k$.
	\end{Rema}
	
	\begin{proof}[Proof of \Cref{thm:stratum_H(2p^q)}]
		As before, we consider a group $G$ of order $n = (2k+1)q$, together with two generators $x,y$ such that $[x, y]$ has order $2k + 1$.
		
		First, observe that $H = \langle [x, y]\rangle \leq G'$ is a cyclic subgroup of $G$ of prime index $q$. We deduce that either $G' = H$ or $G = G'$. The latter case cannot occur. Indeed, it would mean that $G$ is nonsolvable, but, since $|G|$ is odd, the Feit--Thompson theorem (\Cref{thm:feit-thompson}) implies that $G$ is solvable.
		
		We deduce that $G' = H$. In particular, since the commutator subgroup is always normal, we get that $H$ is normal in $G$.
		
		Now, we use \Cref{cor:p_Sylow_order_p_cyclic} to obtain that:
		\[
		G \simeq \ZZ/(2k + 1) \rtimes_d \ZZ/q\ZZ
		\]
		for some $d \in \{1,\cdots, 2k\}$. We know that such a semidirect product exists if and only if $d^q \equiv 1 \mod 2k+1$. Moreover, its commutator subgroup is $\ZZ/(2k+1)\ZZ$ (see \Cref{lem:commutator_metacyclic}): this is possible if and only if $d-1$ is coprime with $2k+1$. In this case, the generators $x = (1, 0)$ and $(y, 0)$ satisfy $[x, y] = (1 - d, 0)$, which has order $2k + 1$. This proves the first part of the statement.
		
		Finally, the existence of $d \in \{2, \dots, 2k\}$ meeting the requirements is equivalent to all prime factors of $2k+1$ being congruent to $1$ modulo $q$ by  \Cref{prop:derived_sgp_semi_direct_product_vfacile}.
	\end{proof}
	Using \Cref{lem:extension_cyclic}, this construction can be bootstrapped to other strata:
	
	\begin{Cor}\label{cor:p^2l+1}
		Let $k \geq 1$ be an integer and let $q$ be an odd prime such that every prime factor of $2k + 1$ is congruent to $1$ modulo $q$. Then, there exist regular origamis in the stratum $\calH(2k^{\ell q})$. In particular, $g = \ell k q + 1 \notin \genera(\infty)$.
	\end{Cor}
	
	\begin{proof}
		The previous proposition shows the existence of a regular origami in the stratum $\calH(2k^q)$ with translation group $\ZZ/(2k + 1)\ZZ \rtimes \ZZ/q\ZZ$.
		
		Now, the hypothesis implies $\gcd(2k + 1, q) = 1$. Thus, \Cref{lem:extension_cyclic} directly shows that regular origamis with translation group $(\ZZ/(2k + 1)\ZZ \rtimes \ZZ/q\ZZ) \times \ZZ/\ell\ZZ$ exist in the stratum $\calH(2k^{\ell q})$.
	\end{proof}
	
	An immediate consequence of this proposition is that, if $\calH(2k^q)$ contains regular origamis, then $2k+1 \equiv 1 \mod q$ and, hence, $q \divides 2k$. In particular, since $q\neq 2$, we must have $q \divides k$. If $k=p$ is a prime number, we get:
	\begin{Cor} \label{cor:Sophie_Germain}
		If $p,q$ are odd primes, there exist regular origamis in the stratum $\mathcal{H}(2p^q)$ if and only if $p = q$ and $2p + 1$ is prime.
	\end{Cor}
	
	\begin{proof}
		Taking $\ell = 1$ in the previous corollary shows existence.
		
		Now, assume a regular origami exists in $\calH(2p^q)$. By \Cref{thm:stratum_H(2p^q)} and the discussion above, we get that $q \divides p$ and, since $p$ is prime, that $p = q$. Now, \Cref{thm:stratum_H(2p^q)} also shows that every prime factor $r$ of $2p + 1$ is congruent to $1$ modulo $p$. This can only happen if $r = 2p + 1$ and $2p + 1$ is itself prime: otherwise we get $r < p$, so $r$ is not congruent to $1$ modulo $p$.
	\end{proof}
	
	\subsection{The stratum \texorpdfstring{$\mathcal{H}(k^{2q})$}{H(k\textasciicircum(2q))}.} \label{s:H(p^2q)}
	We will now study regular origamis in the statum $\calH(k^{2q})$, where $q$ is an odd prime, and $k$ is odd. 
	Our goal is to show:
	
	\begin{Thm}\label{thm:stratum_H(k^2q)}
		Let $k$ be an odd integer and $q$ be an odd prime. Then, regular origamis exist in the stratum $\calH(k^{2q})$ if and only if $q = 3$ and either:
		\begin{itemize}[wide, labelindent=0pt]
			\item $k \equiv 1 \mod{4}$ and every prime factor of $(k + 1) / 2$ is congruent to $1$ modulo $3$; or
			\item $k \equiv 3 \mod{4}$ and every prime factor of $(k + 1) / 4$ is congruent to $1$ modulo $3$.
		\end{itemize}
		Moreover, any regular origami belonging to such strata has a translation group isomorphic to either $(\ZZ/((k + 1)/2)\ZZ \times \ZZ/2\ZZ \times \ZZ/2\ZZ) \rtimes \ZZ/3\ZZ$ (in the former case) or $(\ZZ/((k + 1)/4)\ZZ \times \Q_8) \rtimes \ZZ/3\ZZ$ (in the latter case).
	\end{Thm}
	
	\begin{Rema}\label{rema:origami_with_dihedral}
		When either $q=2$ or $k$ is even, we already know that there exist regular origamis in $\calH(k^{2q})$, constructed via semidirect products, see \Cref{sec:examples_origamis}.
	\end{Rema}
	
	An equivalent formulation of \Cref{thm:stratum_H(k^2q)} is the following:
	
	\begin{Thm}\label{thm:stratum_H(k^2q)_group}
		Let $k$ be an odd integer and let $q$ be an odd prime. Then, a group $G$ of order $n = 2(k + 1)q$, together with two generators $x, y \in G$ such that $[x, y]$ has order $k + 1$, exists if and only if $q = 3$ and either:
		\begin{itemize}[wide, labelindent=0pt]
			\item $k \equiv 1 \mod{4}$ and every prime factor of $(k + 1) / 2$ is congruent to $1$ modulo $3$; or
			\item $k \equiv 3 \mod{4}$ and every prime factor of $(k + 1) / 4$ is congruent to $1$ modulo $3$.
		\end{itemize}
		Moreover, $G$ is isomorphic to either $(\ZZ/((k + 1)/2)\ZZ \times \ZZ/2\ZZ \times \ZZ/2\ZZ) \rtimes \ZZ/3\ZZ$ (in the former case) or $(\ZZ/((k + 1)/4)\ZZ \times \Q_8) \rtimes \ZZ/3\ZZ$ (in the latter case).
	\end{Thm}
	
	As before, we will use the notation $H = \langle [x, y]\rangle$. We start by showing that $H \neq G'$, by contradiction.
	
	\begin{Lem}\label{lem:H_not_G'}
		In the context of \Cref{thm:stratum_H(k^2q)_group}, we have $H < G'$.
	\end{Lem}
	
	\begin{proof}
		Assume by contradiction that $H = G'$. Observe that $G / G'$ is abelian and has order $2q$, so it is isomorphic to $\ZZ/2q\ZZ$. Thus, $G$ is a metacyclic group with presentation:
		\[
		\langle a, b \st a^{k+1} = 1, b^{2q} = a^j, [b, a] = a^{d-1}\rangle,
		\]
		where $d^{2q} \equiv 1 \mod{k+1}$ and $(k+1) \divides j(d - 1)$.  The first condition implies that $d$ is coprime with $k+1$. Since $k+1$ is even, we get that $r$ is odd. As a consequence, $d-1$ is even and the commutator subgroup of $G$, which is generated by $a^{d-1}$, is strictly contained in $H$. This is a contradiction.
	\end{proof}
	
	Since $(G : H) = 2q$, we deduce that $(G : G')$ divides $2q$, but does not equal $2q$, so this quantity is either $1$, $2$ or $q$. We will analyze each of these cases separately, obtaining that the only admissible case is $(G : G') = q$. Unlike the case of $\calH(2k^q)$, the subgroup $H$ is \emph{never} normal in $G$ (see \Cref{lem:stratum_H(k^2q)_index_q_not_normal}). Thus, we need to resort to more intricate devices. 
	
	\subsubsection{$G$ is perfect} We start by assuming that $G$ is perfect, that is, $G' = G$ and we will derive a contradiction.
	
	Throughout this section, we will not assume that $k$ is odd, so we will prove a more general statement than needed for \Cref{thm:stratum_H(k^2q)_group}. That is, our goal is to show:
	
	\begin{Prop} \label{prop:stratum_H(k^2q)_perfect}
		Let $q$ be an odd prime. Then, a perfect group $G$, together with two generators $x, y \in G$, such that $H = \langle[x, y]\rangle$ has index $2q$, does not exist.
	\end{Prop}
	
	\begin{Rema}
		Since there exist $2$-generated perfect groups, there exist regular origamis in some strata whose translation group is perfect. A notable example is the case of $\PSL(2, p)$, treated in \Cref{sec:PSL}.
	\end{Rema}
	
	We will assume the existence of such a group $G$. In the spirit of \Cref{thm:stratum_H(k^2q)_group}, we write $|G| = 2(k + 1)q$. Observe that this quantity has at least \emph{three} distinct prime factors: $2$, $q$, and at least one more prime $p$. Indeed, if this were not the case, $G$ would be solvable by Burnside's theorem (\Cref{thm:Burnside-p^a_q^b}). This is impossible as nontrivial perfect groups are nonsolvable. In particular, we have $p \divides (k + 1)$. We will use this fact several times throughout the proof of \Cref{prop:stratum_H(k^2q)_perfect}.
	In fact, we can also show that the order of $\Z(G)$ is not a multiple of $p$. More precisely:
	
	\begin{Lem}\label{lm:center}
		In the context of \Cref{thm:stratum_H(k^2q)_group}, we have that $|\Z(G)|$ divides $2q$.
	\end{Lem}
	
	The proof of this lemma relies on the following result:
	
	\begin{Thm}[{\cite[Corollary 5.9]{Isaacs_Gp_theory}}]\label{thm:Cor5.9_Isaacs}
		Suppose that $G$ is a finite group and that $(G : \Z(G)) = \ell$. Then, the $\ell$-th power of every commutator is the identity.
	\end{Thm}
	
	\begin{proof}[Proof of~\Cref{lm:center}]
		Let $\ell = (G : \Z(G))$. By \Cref{thm:Cor5.9_Isaacs}, the $\ell$-th power of $[x, y]$ is trivial. Since the order of $[x, y]$ is $k + 1$, $\ell$ must be a multiple of $k + 1$, say $\ell =  t(k + 1)$. Hence,
		\[
		|\Z(G)| = \frac{|G|}{(G : \Z(G))} = \frac{2q}{t}. \qedhere
		\]
	\end{proof}
	
	In fact, we will show that up to taking the quotient by $\Z(G)$, we can assume $\Z(G) = 1$. This is a consequence of:
	\begin{Lem}\label{lem:center_in_H} In the context of \Cref{prop:stratum_H(k^2q)_perfect}, we have
		$\Z(G) \leq H$.
	\end{Lem}
	
	Indeed, if we assume this Lemma to be true, then the group $G/\Z(G)$:
	\begin{itemize}
		\item is generated by $x \Z(G)$ and $y \Z(G)$;
		\item is perfect, as it is a quotient of a perfect group;
		\item is centerless, from Grün's Lemma (\Cref{lem:Grun's}); and
		\item admits a cyclic subgroup $\langle [x \Z(G), y \Z(G)]\rangle = H/\Z(G)$ of index $2q$. 
	\end{itemize}
	Thus, if a group $G$ as in \Cref{prop:stratum_H(k^2q)_perfect} exists, then a centerless group $G/\Z(G)$ as in \Cref{prop:stratum_H(k^2q)_perfect} also exists. Hence, we can assume that $G$ is centerless.
	
	Before proving \Cref{lem:center_in_H}, we need two preliminary results.
	
	\begin{Thm}[{\cite[Theorem 5.18]{Isaacs_Gp_theory}}]\label{thm:Thm5.18_Isaacs}
		Let $P$ be an abelian Sylow $p$-subgroup of a finite group $G$. Then,
		\[ G' \cap P \cap \Z(\N_G(P)) = \{1 \}.\]
	\end{Thm}
	
	\begin{Lem} \label{lem:index_N_G(H)_in_G_and_H_in_N_G(H)}
		In the context of \Cref{prop:stratum_H(k^2q)_perfect}, fix a prime number $p$ different from $2$ and $q$ dividing $k + 1$. Then, there exists a unique Sylow $p$-subgroup of $G$ contained in $H$. Moreover, we have $\N_G(P) = \N_G(H)$, and
		\[
		(G : \N_G(H)) = q \quad \text{and} \quad (\N_G(H) : H) = 2.
		\]
		In particular, $q \equiv 1 \mod p$, so $q > 3$.
	\end{Lem}
	\begin{proof}
		Let $p$ be a prime number $p$ different from $2$ and $q$ dividing $k + 1$. There exists a Sylow $p$-subgroup $P$ of $H$, which is unique as $H$ is cyclic. Since $(G : H) = 2q$ and $p$ are coprime, $P$ is also a Sylow $p$-subgroup of $G$. Furthermore, since $P \leq H$ and $H$ is abelian, every element of $H$ normalizes $P$, so we have $H \leq \N_G(P) \leq G$. We will show that $(\N_G(P) : H) = 2$ and that $\N_G(H) = \N_G(P)$. We start by showing that both inclusions are strict.
		
		\smallbreak
		
		If $\N_G(P) = G$, then $P$ is normal in $G$. By~\Cref{lm:action}, 
		we have $P \leq \Z(G)$.
		This is impossible since $p$ does not divide $|\Z(G)|$, as $|\Z(G)|$ divides $2q$ by \Cref{lm:center}.
		
		\smallbreak
		
		If $\N_G(P) = H$, we have $\Z(\N_G(P)) = \Z(H) = H$ since $H$ is abelian. Moreover, we have $G' \cap P \cap \Z(\N_G(P)) = P$ since $G$ is perfect. This contradicts~\Cref{thm:Thm5.18_Isaacs}.
		
		\smallbreak
		
		We deduce that $\N_G(P)$ is an index-$q$ subgroup of $G$ containing $H$. Indeed, the number of Sylow $p$-subgroups of $G$ is exactly $(G : \N_G(P))$, and this number is congruent to $1$ modulo $p$ from \Cref{thm:Sylow}. Since $(G : H) = 2q$ and $H < \N_G(P)$, we obtain that $(G : \N_G(P))$ is either $2$ or $q$. Since $2 \not\equiv 1 \mod{p}$, the only possibility is that $(G : \N_G(P)) = q$ and that $q \equiv 1 \mod{p}$. Hence, $(\N_G(P) : H) = 2$ and $q > p > 2$, so $q > 3$.
		
		Finally, we show that $\N_G(H)= \N_G(P)$. Indeed, since $(\N_G(P) : H) = 2$, we have that $H \triangleleft \N_G(P)$ by \Cref{lem:index_two_normal}. Thus, $\N_G(P) \leq \N_G(H)$. Conversely, we have $\N_G(H) \leq \N_G(P)$ since $P \leq H$ is a characteristic subgroup of $H$ and, therefore, conjugation by an element $g \in \N_G(H)$, which induces an automorphism of $H$, also normalizes $P$. This proves that $\N_G(H) = \N_G(P)$.
	\end{proof}
	
	We can now prove \Cref{lem:center_in_H}.
	\begin{proof}[Proof of \Cref{lem:center_in_H}]
		Let $p$ be a prime number different from $2$ and $3$, and let $P$ be a Sylow $p$-subgroup of $G$ contained in $H$. We have $\Z(G) \leq \C_G(P) \leq \N_G(P)$.
		
		Now, assume by contradiction that $\Z(G) \not\leq H$ and take $t \in \Z(G) \setminus H$. Since $t \in \N_G(P)$ and $(\N_G(P):H) = 2$ by \Cref{lem:index_N_G(H)_in_G_and_H_in_N_G(H)}, we have $\N_G(P) = \langle H,t \rangle$. Moreover, $t \in \C_G(P)$ and also $H \leq \C_G(P)$ since $P \leq H$ and $H$ is abelian, so we deduce that
		\[
		\N_G(P) = \langle H, t \rangle \leq \C_G(P) \leq \N_G(P).
		\]
		Hence, $P \leq \Z(\C_G(P)) = \Z(\N_G(P))$, which again contradicts \Cref{thm:Thm5.18_Isaacs}. We conclude that $\Z(G) \leq H$.
	\end{proof}
	
	As previously discussed, we may now assume that $\Z(G) = 1$. We will investigate the number of elements of order $q$ to prove that $q \leq 3$, which will contradict $q > 3$ from \Cref{lem:index_N_G(H)_in_G_and_H_in_N_G(H)}. We first show:
	
	\begin{Lem}\label{lem:intersection_H_tHt-1}
		In the context of \Cref{prop:stratum_H(k^2q)_perfect}, and assuming that $G$ is centerless, we have $H \cap t H t^{-1} = \{1\}$ for any $t \in G \setminus \N_{G}(H)$.
	\end{Lem}
	\begin{proof}
		We first prove that $K = H \cap t H t^{-1}$ is normal by showing that $\N_G(K) = G$. Since $(G : \N_G(H)) = q$ is prime, we have $G = \langle \N_G(H), t\rangle$, so it is enough to establish that $\N_G(H) \leq \N_G(K)$ and that $t \in \N_G(K)$.
		
		Let $s \in \N_G(H)$. We have that $K \leq H$ and that $s K s^{-1} \leq sHs^{-1} = H$. Since $s K s^{-1}$ and $K$ are subgroups of the same order of the cyclic subgroup $H$, we deduce that $s K s^{-1}= K$ by \Cref{thm:fundamental_cyclic_groups}. Hence, $s \in \N_G(K)$.
		
		Now, $K \leq t H t^{-1}$ and $t K t^{-1} \leq t H t^{-1}$. Again, $t K t^{-1}$ and $K$ are subgroups of the same order of the cyclic subgroup $t H t^{-1}$, so $t K t^{-1} = K$ by \Cref{thm:fundamental_cyclic_groups}. We get that $t \in \N_G(K)$.
		
		Finally, we have that $tHt^{-1} \cap H$ is a cyclic normal subgroup of $G$, and therefore it is central in $G$ by
		\Cref{lm:action}. Since $\Z(G) = 1$, we get that $tHt^{-1} \cap H$ is trivial.
	\end{proof}
	
	The following is a somewhat direct consequence of \Cref{lem:intersection_H_tHt-1}.
	
	\begin{Lem} \label{lem:intersection_single_element}
		In the context of \Cref{prop:stratum_H(k^2q)_perfect}, and assuming that $G$ is centerless, define $H_2 = \{h \in H \st h^2 = 1\}$, let $t \in G \setminus \N_G(H)$, and take $0 \leq \beta < \alpha$ with $t^{\alpha - \beta} \notin \N_G(H)$. Then,
		\[
		t^{\alpha}(\N_G(H) \setminus H_2)t^{-\alpha} \cap t^{\beta}(\N_G(H) \setminus H_2)t^{-\beta}
		\]
		contains at most a single element.
	\end{Lem}
	
	\begin{proof}
		Recall that $(\N_G(H) : H) = 2$ from \Cref{lem:index_N_G(H)_in_G_and_H_in_N_G(H)}.
		
		Let $u = t^{\alpha - \beta}$ and $K = \N_G(H) \cap u \N_G(H) u^{-1}$. If $s \in K$, we have $s^2 = 1$. Indeed, $s^2 \in H$ by \Cref{lem:index_two_order_two}. Similarly, $s^2 \in u H u^{-1}$. Hence, $s^2 = 1$ by \Cref{lem:intersection_H_tHt-1}.
		
		Next, we will show that the set
		\[
		S = (\N_G(H) \setminus H) \cap u (\N_G(H) \setminus H) u^{-1} \subseteq K
		\]
		has at most a single element. If $s, s' \in S$, we have that $s s' \in H \cap u H u^{-1}$ by \Cref{lem:index_two_order_two}, so $ss' = 1$ by \Cref{lem:intersection_H_tHt-1}. Hence, $s' = s^{-1} = s$. We get that $|S| \leq 1$.
		
		Now, let
		\[
		T = (\N_G(H) \setminus H_2) \cap u (\N_G(H) \setminus H_2) u^{-1} \subseteq K.
		\]
		We will show that $T \subseteq S$, so $T$ also has at most a single element. Let $s \in T$.
		
		By definition, $s \notin H_2$, so $s \notin H$ as every element of $K$ has order at most $2$. Similarly, $s \notin u H_2 u^{-1}$, so $s \notin u H u^{-1}$. We deduce that $s \in S$.
		
		Finally, we have that
		\[
		t^{\alpha}(\N_G(H) \setminus H_2)t^{-\alpha} \cap t^{\beta}(\N_G(H) \setminus H_2)t^{-\beta} = t^{\beta} T t^{-\beta},
		\]
		so the latter set also has at most a single element.
	\end{proof}
	
	\Cref{lem:intersection_H_tHt-1} also allows us to estimate the number of elements of order $q$ in $G$.
	\begin{Lem} \label{lem:order_q}
		In the context of \Cref{prop:stratum_H(k^2q)_perfect}, and assuming that $G$ is centerless, we have that $q^2$ does not divide $|G|$. Moreover, the number of Sylow $q$-subgroups of $G$ is at least $q + 1$, and each of them is isomorphic to $\ZZ/q\ZZ$. In particular, the number of order-$q$ elements of $G$ is at least $q^2 - 1$.
	\end{Lem}
	\begin{proof}
		Let $t \in G \setminus \N_{G}(H)$. Since $t Ht^{-1} \cap H = \{1\}$ by \Cref{lem:intersection_H_tHt-1} the cardinality of the set $(t H t^{-1}) H$ is $|H|^2$. Indeed, by \Cref{lem:order_HK}:
		\[
		\frac{|H| |t H t^{-1}|}{|H \cap t H t^{-1}|} = |H|^2.
		\]
		Hence, $|H|^2 \leq |G|$. Dividing by $|H|$ yields
		\[ 
		k + 1 = |H| \leq (G : H) = 2q.
		\]
		Since $|G| = 2(k + 1) q$, we obtain that $q^2$ does not divide $|G|$. Indeed, we know that $k + 1$ has a prime factor $p$ different from $2$ and $q$, and
		\[
		\frac{k + 1}{p} \leq \frac{2 q}{p} < q.
		\]
		Therefore, any Sylow $q$-subgroup $Q$ of $G$ must be isomorphic to $\ZZ/q\ZZ$.
		
		Now, the number of Sylow $q$-subgroups is congruent to $1$ modulo $q$. However, this number cannot be $1$: this would mean that $Q$ is normal and, hence, central by \Cref{lm:action}. This is impossible as $\Z(G) = 1$. In particular, the number of Sylow $q$-subgroups of $G$ is at least $q + 1$.
		
		Finally, since every nontrivial element of $\ZZ/q\ZZ$ is generating, the intersection of two distinct Sylow $q$-subgroups is trivial. Thus, the number of order-$q$ elements of $G$ is at least $(q + 1)(q - 1) = q^2 - 1$.
	\end{proof}
	
	We can now finish the proof of \Cref{prop:stratum_H(k^2q)_perfect}.
	
	\begin{proof}[Proof of \Cref{prop:stratum_H(k^2q)_perfect}]
		First, recall that we can assume that $G$ is centerless. Fix an order-$q$ element $t \in G$.
		
		We have $t \notin \N_G(H)$ since $q \ndivides 2(k + 1) = |\N_G(H)|$ as $q^2 \ndivides 2(k + 1)q = |G|$ by \Cref{lem:order_q}. Similarly, $t^\alpha \notin \N_G(H)$ for every $1 \leq \alpha \leq q - 1$, as the order of this element is also $q$.
		
		Consider the set $S = \bigcup_{\alpha = 0}^{q - 1} t^\alpha \N_G(H) t^{-\alpha}$. We claim that
		\[
		|S| \geq 1 + |G| - \frac{q(q + 3)}{2}.
		\]
		Indeed, set $H_2 = \{h \in H \st h^2 = 1\}$. From \Cref{thm:fundamental_cyclic_groups}, we have $|H_2| \leq 2$. Now, we estimate the elements of $S$ by excising $H_2$ from $\N_G(H)$ and using \Cref{lem:intersection_single_element}:
		\begin{align*}
			|S| &\geq 1 + \sum_{\alpha = 0}^{q - 1} \left(|t^{\alpha} (\N_G(H) \setminus H_2) t^{-\alpha}|\vphantom{\sum_{\alpha=0}^{q-1}}\right. \\
			&\left.\quad -\sum_{\beta = 0}^{\alpha - 1} |t^{\alpha} (\N_G(H) \setminus H_2) t^{-\alpha} \cap t^{\beta} (\N_G(H) \setminus H_2) t^{-\beta}|\right) \\
			&\geq 1 + \sum_{\alpha = 0}^{q - 1} (|\N_G(H)| - 2 - \alpha) = 1 + q|\N_G(H)| - \frac{q(q + 3)}{2} = 1 + |G| - \frac{q(q + 3)}{2},
		\end{align*}
		where the ``$1$'' corresponds to counting the trivial element, and where we used that $(G : \N_G(H)) = q$ from \Cref{lem:index_N_G(H)_in_G_and_H_in_N_G(H)}.
		
		Finally, \Cref{lem:order_q} shows that $G$ contains at least $q^2 - 1$ elements of order $q$. Since no element of $S$ has order $q$ (using that $q$ does not divide $|\N_G(H)|$), we get:
		\[
		(q^2 - 1) + \left(1 + |G| - \frac{q(q + 3)}{2}\right) \leq |G|.
		\]
		Hence,
		\[ 
		q^2 \leq \frac{q(q + 3)}{2}.
		\]
		This inequality holds if and only if $0 \leq q \leq 3$, which contradicts that $q > 3$ from \Cref{lem:index_N_G(H)_in_G_and_H_in_N_G(H)}.
	\end{proof}
	
	\subsubsection{\texorpdfstring{$G'$}{G'} has index two} We now assume that $G'$ as index two inside $G$, and we will again derive a contradiction. That is, our goal is to show:
	
	\begin{Lem} \label{prop:stratum_H(k^2q)_index_two}
		Let $k$ be an odd integer and let $q$ be an odd prime. Then, a group $G$ of order $n = 2(k + 1)q$ with $(G : G') = 2$, together with two generators $x, y \in G$ such that $[x, y]$ has order $k + 1$, does not exist.
	\end{Lem}
	
	We start by showing that such a group must admit a normal $2'$-Hall subgroup:
	\begin{Lem} \label{lem:stratum_H(k^2q)_index_two_2'_Hall}
		In the context of \Cref{prop:stratum_H(k^2q)_index_two}, $G$ contains a normal $2'$-Hall subgroup.
	\end{Lem}
	
	\begin{proof}
		We will use the Frobenius $p$-complement theorem (\Cref{thm:Frobenius_p_complement}), that is, we need to show that, for every $2$-subgroup $X \leq G$, we have that $\N_G(X) / \C_G(X)$ is also a $2$-group.
		
		Recall that $H = \langle [x, y]\rangle$ and write $k + 1 = 2^\alpha\ell$, with $\ell$ odd. By \Cref{thm:fundamental_cyclic_groups}, $H$ admits a cyclic subgroup $L \leq H$ of order $2^\alpha$.
		
		Observe that $(G : L) = 2\ell q$ is even, but not divisible by $4$. Thus, \Cref{lem:cyclic_index_two} shows that $X$ admits a cyclic subgroup of index two. Hence, $X$ is isomorphic to one of the groups in \Cref{thm:Burnside}.
		
		Recall that $\N_G(X) / \C_G(X)$ injects into $\Aut(X)$ (\Cref{rk:normalizer_centralizer_automorphism}). By \Cref{lem:Burnside_Aut}, we have that $\Aut(X)$ is $2$-group, except when $X \simeq \ZZ/2\ZZ \times \ZZ/2\ZZ$ and when $X \simeq \Q_8$. Thus, we only need to focus on these two cases.
		
		Assume then that $X \simeq \ZZ/2\ZZ \times \ZZ/2\ZZ$ or $X \simeq \Q_8$. Since $X$ is not cyclic, it is not contained in $G'$. Indeed, $L$ is cyclic and is a Sylow $2$-subgroup of $G'$, so every $2$-subgroup of $G'$ is cyclic by Sylow's theorems (\Cref{thm:Sylow}) and \Cref{thm:fundamental_cyclic_groups}. Then, $(G : G') = 2$ implies that $G = XG'$. Defining $K = X \cap G'$ and using that $G' \triangleleft G$, the second isomorphism theorem shows that
		\[
		(X : K) = (X G' : G') = (G : G') = 2.
		\]
		Since every index-two subgroup of $X$ is cyclic, we also obtain that $K$ is cyclic.
		
		Let $s \in K$ be a generator of $K$. Let $t \in X \setminus K = X \setminus G'$. Since $(X : K) = 2$ is prime, we have that $X = \langle K, t\rangle$, so $X = \langle s, t\rangle$.
		
		Now, let $\varphi \in \N_G(X) / \C_G(X)$ seen a subgroup of $\Aut(X)$. We claim that $\varphi$ preserves $K$ and $X \setminus K$. To see this, observe that $\varphi$ acts on $X$ by conjugation by an element of $G$, so there exists an inner automorphism $\overline{\varphi} \in \Inn(G)$ such that $\overline{\varphi}|_X = \varphi$. Since $G'$ is normal in $G$, $\overline{\varphi}$ preserves $G'$, so $\varphi$ preserves $K$. The second equality follows from the first, as $\varphi$ is bijective.
		
		Finally, we check the two particular cases:
		
		\begin{enumerate}[wide, labelindent=0pt, label=\underline{Case \arabic{*}:}, itemsep=1ex]
			\item When $X \simeq \ZZ/2\ZZ \times \ZZ/2\ZZ$, the group $\N_G(X)/\C_G(X)$ has at most $2$ elements. Indeed, $\varphi(s) = s$ and $\varphi(t) \in \{t, ts\}$.
			\item When $X \simeq \Q_8$, the order of $\N_G(X)/\C_G(X)$ divides $8$. Indeed, $s$ has order $4$, so $\varphi(s) \in \{s,s^3\}$. Moreover, $\varphi(t) \in \Q_8 \setminus \langle s \rangle = \{t, st,s^2t, s^3t\}$. We get $\N_G(X)/\C_G(X) \leq \ZZ/2\ZZ \times \ZZ/4\ZZ$.
		\end{enumerate}
		
		From \Cref{thm:Frobenius_p_complement}, we conclude that $G$ has a normal $2$-complement.
	\end{proof}
	
	We can now finish the proof of \Cref{prop:stratum_H(k^2q)_index_two}:
	\begin{proof}[Proof of \Cref{prop:stratum_H(k^2q)_index_two}]
		By \Cref{lem:stratum_H(k^2q)_index_two_2'_Hall}, $G$ contains a normal $2'$-Hall subgroup $N$. Since $2^2 \ndivides 2 = (G : G')$, \Cref{lem:p_Sylow_order_p_Hall} shows that
		\[
		G \simeq N \rtimes \ZZ/2\ZZ,
		\]
		where $|N|$ is odd. In particular, $|G|$ is not divisible by $4$, which contradicts that $|G| = 2(k + 1)q$ for odd $k$.
	\end{proof}
	
	\subsubsection{$G'$ has index $q$} We finally focus on the case where $G'$ has index $q$ inside $G$. We will show that regular origamis can exist in the stratum $\calH(k^{2q})$ under very special conditions:
	\begin{Prop}\label{prop:stratum_H(k^2q)_index_q}
		Let $k$ be an odd integer and let $q$ be an odd prime. Then, a group $G$ of order $n = 2(k + 1)q$ with $(G : G') = q$, together with two generators $x, y \in G$ such that $[x, y]$ has order $k + 1$, exists if and only if $q = 3$ and either:
		\begin{itemize}[wide, labelindent=0pt]
			\item $k \equiv 1 \mod{4}$ and every prime factor of $(k + 1) / 2$ is congruent to $1$ modulo $3$; or
			\item $k \equiv 3 \mod{4}$ and every prime factor of $(k + 1) / 4$ is congruent to $1$ modulo $3$.
		\end{itemize}
		Moreover, $G$ is isomorphic to either $(\ZZ/((k + 1)/2)\ZZ \times \ZZ/2\ZZ \times \ZZ/2\ZZ) \rtimes \ZZ/3\ZZ$ (in the former case) or $(\ZZ/((k + 1)/4)\ZZ \times \Q_8) \rtimes \ZZ/3\ZZ$ (in the latter case).
	\end{Prop}
	
	\begin{Rema}
		Any regular origami in $\calH(k^{2q})$, for odd $k$ and odd prime $q$, has genus $g$ satisfying $g - 1 = qk = 3k$, so $3 \divides (g - 1)$. Since the genera covered by the statement of \Cref{thm:no_regular_origamis} never satisfy this property, the existence of such origamis is not an obstruction for this theorem.
	\end{Rema}

	We start by remarking that $(G' : H) = 2$, so $H$ is normal in $G'$ by \Cref{lem:index_two_normal}. We will use this fact several times. Nevertheless, $H$ is not normal in $G$:
	
	\begin{Lem}\label{lem:stratum_H(k^2q)_index_q_not_normal}
		Let $G$ be a group as in \Cref{lem:stratum_H(k^2q)_structure_G'}. Then, $H = \langle[x, y]\rangle$ is not normal in $G$.
	\end{Lem}
	
	\begin{proof}
		Assume that $H$ is normal. Consider the short exact sequence
		\[
		1 \to G' \to G \to \ZZ/q\ZZ \to 1.
		\]
		Since $(G' : H) = 2$, the quotient by $H$ yields the short exact sequence:
		\[
		1 \to \ZZ/2\ZZ \simeq G'/H \to G/H \to \ZZ/q\ZZ \to 1.
		\]
		We deduce that $G/H$ is a metacyclic group with presentation:
		\[
		G/H = \langle a, b \st a^2 = 1, b^q = a^\varepsilon \text{ and } [b, a] = a^{d-1}\rangle,
		\]
		where $d^q \equiv 1 \mod 2$ and $2 \divides \varepsilon(d - 1)$. In particular, $d$ is odd, so $d - 1$ is even and the group is abelian. Hence, $(G/H)'$ is trivial. We deduce that $G'/H \simeq (G/H)'$ is trivial, which is a contradiction.
	\end{proof}
	
	Now, we can greatly restrict the structure of $G'$:
	
	\begin{Lem}\label{lem:stratum_H(k^2q)_structure_G'}
		In the context of \Cref{prop:stratum_H(k^2q)_index_q}, we have $G' \simeq \ZZ/\lambda\ZZ \times L$, where either $L \simeq \ZZ/2\ZZ \times \ZZ/2\ZZ$ or $L \simeq \Q_8$ and $\lambda$ is odd. Furthermore, $L$ is normal in $G$.
	\end{Lem}
	\begin{proof}
		Write $k+1 = |H| = 2^\alpha \lambda$ (where $\alpha \geq 1$ since $k$ is odd by assumption). Let $K$ be the unique subgroup of $H$ of order $\lambda$. Notice that $K$ is normal in $G'$, since it is a characteristic subgroup of $H$, which has index two in $G'$. As a consequence, and again using that $(G:H)=2$, we obtain that $K$ is a normal $2'$-Hall subgroup of $G'$. In particular, the Schur--Zassenhaus theorem (\Cref{thm:Schur_Zassenhaus}) shows that:
		\[
		G' \simeq \ZZ/\lambda\ZZ \rtimes L,
		\]
		where $L$ is a $2$-group of order $2^{\alpha+1}$.
		In fact, a normal Hall subgroup is characteristic \cite[Exercise 5.31]{Rotman_Groups}, so $K$ is characteristic in $G'$ and normal in $G$. In particular, we can use \Cref{lm:action} to conclude that $K \leq \Z(G')$.
		This implies that the semidirect product is actually a direct product. Moreover, the $2$-group $L$ contains a cyclic subgroup of index two, so it is one of the groups of \Cref{thm:Burnside}. Furthermore, since $\lambda$ and $|L|$ are coprime, we have
		\[
		\Aut(G') \simeq \Aut(\ZZ/\lambda\ZZ \times L) \simeq \Aut(\ZZ/\lambda\ZZ) \times \Aut(L),
		\]
		so, in particular, $L$ is characteristic in $G'$ and, hence, normal in $G$ by \Cref{lem:characteristic_in_normal_is_normal}.
		
		Assume now by contradiction that $L$ is not isomorphic to $\ZZ/2\ZZ \times \ZZ/2\ZZ$ or $\Q_8$. Using \Cref{lem:Burnside_characteristic}, it contains a characteristic subgroup $U$ of index two. Since $L$ is normal in $G$, we deduce that $U$ is normal in $G$ by \Cref{lem:characteristic_in_normal_is_normal}. Thus, the group $G/U$ lies in the short exact sequence:
		\[
		1 \to G'/U \to G/U \to G/G' \to 1
		\]
		Since $G'/U \simeq (\ZZ/\lambda \ZZ \times L)/U \simeq \ZZ/\lambda \ZZ \times \ZZ/2\ZZ \simeq \ZZ/2\lambda \ZZ$ and $G/G' \simeq \ZZ/q\ZZ$ is cyclic, we have:
		\[
		1 \to \ZZ/2\lambda\ZZ \to G/U \to \ZZ/q\ZZ \to 1.
		\]
		Therefore, $G/U$ is a metacyclic group with presentation:
		\[
		\langle x, y \st x^{2\lambda} = 1, y^q = x^r \text{ and } [y,x] = x^{d-1} \rangle,
		\]
		where $d^q = 1 \mod{2\lambda}$. In particular, $d$ is odd. Consequently, $d-1$ is even and $2 \divides \gcd(d-1, 2\lambda)$. In particular, from \Cref{lem:commutator_metacyclic}, we deduce that $(G/U)'$ is strictly contained in $G'/U \simeq \ZZ/2\lambda\ZZ$. This is a contradiction since $(G/U)' \simeq G'/U$.
	\end{proof}
	
	According to \Cref{lem:stratum_H(k^2q)_structure_G'}, we only have to distinguish two cases for $G'$. When $G' \simeq \ZZ/\lambda \ZZ/2\ZZ \times \ZZ/2\ZZ$, we have that $|G| = 2(k + 1)q = (4\lambda) \cdot q$, so $\lambda = (k + 1)/2$. This number is an odd integer, so we get that $k \equiv 1 \mod{4}$. Similarly, when $G' \simeq \ZZ/\lambda\ZZ \times \Q_8$, we have that $|G| = 2(k + 1)q = (8\lambda) \cdot q$, so $\lambda = (k + 1)/4$. This number is an odd integer, so we get $k \equiv 3 \mod{8}$ (and, in particular, $k \equiv 3 \mod{4}$).
	
	We will now show the desired conditions on $q$ and $\lambda$, and group structure for $G$:
	
	\begin{Lem}\label{lem:quaternion_klein_x_cyclic}
		Let $G$ be a group as in \Cref{prop:stratum_H(k^2q)_index_q}. Then, $q=3$ and every prime factor of $\lambda$ is congruent to $1$ modulo $3$. Furthermore, $G \simeq (\ZZ/\lambda\ZZ \times L) \rtimes \ZZ/3\ZZ$, where $L \simeq \ZZ/2\ZZ \times \ZZ/2\ZZ$ or $L \simeq \Q_8$.
	\end{Lem}
	
	\begin{proof}
		Write $G' = \ZZ/\lambda\ZZ \times L$, with $L \simeq \ZZ/2\ZZ \times \ZZ/2\ZZ$ or $L \simeq \Q_8$, as per \Cref{lem:stratum_H(k^2q)_structure_G'}.
		
		We start by proving that $q = 3$. On the one hand, observe that, since $H$ is not normal in $G$ by \Cref{lem:stratum_H(k^2q)_index_q_not_normal}, we have $\N_G(H) < G$. Moreover, since $(G' : H) = 2$, $H$ is normal in $G'$ by \Cref{lem:index_two_normal}, and $G' \leq \N_G(H)$. We deduce $G' \leq \N_G(H) < G$. Since $(G : G') = q$ is prime, we get that $\N_G(H) = G'$. Thus,
		\[
		(G : \N_G(H)) = (G: G') = q,
		\]
		so the number of conjugacy classes of $H$ inside $G$ is exactly $q$.
		
		On the other hand, the number of cyclic subgroups of index two inside $L$ is exactly three: either $\langle (1, 0)\rangle$, $\langle (0, 1)\rangle$, and $\langle (1, 1)\rangle$ if $L \simeq \ZZ/2\ZZ \times \ZZ/2\ZZ$; or $\langle\qi\rangle$, $\langle\qj\rangle$, and $\langle\qk\rangle$ if $L \simeq \Q_8$. This also holds inside $G' \simeq \ZZ/\lambda\ZZ \times L$, since $\lambda$ is odd. Therefore, the number of conjugacy classes of $H$ inside $G'$ is at most $3$. We deduce that $q \leq 3$, so $q = 3$ as it is an odd prime.
		
		Now, $L$ is normal $G$ by \Cref{lem:stratum_H(k^2q)_structure_G'}. The quotient $K = G/L$ is part of the short exact sequence:
		\[
		1 \to \ZZ/\ZZ\lambda \to K \to \ZZ/3\ZZ \to 1,
		\]
		so $K$ is a metacyclic group with presentation
		\[
		K = \langle a, b \st a^\lambda = 1, b^3 = a^i \text{ and } [b, a] = a^{d - 1}\rangle,
		\]
		where $d^3 \equiv 1 \mod \lambda$ and $\lambda \divides i(d - 1)$. Since $H = \langle[x, y]\rangle$ contains $\ZZ/\lambda\ZZ$, we deduce that $[xL, yL] \in K'$ generates the group $\ZZ/\lambda\ZZ \leq K$. Moreover, we have $K' = \langle a^{d - 1}\rangle$ by \Cref{lem:commutator_metacyclic}, so we deduce that $\gcd(d - 1, \lambda) = 1$. This gives the desired conditions on $\lambda$ by \Cref{prop:derived_sgp_semi_direct_product_vfacile}. 
		
		Finally, since $G'$ has order $2\lambda$ or $4\lambda$, with $\lambda \equiv 1 \mod{3}$, and $(G : G') = q = 3$, the Schur--Zassenhaus theorem (\Cref{thm:Schur_Zassenhaus}) implies that
		\[
		G \simeq G' \rtimes \ZZ/3\ZZ \simeq (\ZZ/\lambda\ZZ \times L) \rtimes \ZZ/3\ZZ,
		\]
		as desired.
	\end{proof}
	
	To complete the proof of \Cref{prop:stratum_H(k^2q)_index_q}, we exhibit groups, together with two generators, satisfying the required properties:
	
	\begin{proof}[Proof of \Cref{prop:stratum_H(k^2q)_index_q}]
		Let $k \geq 1$ be odd and assume that it satisfies the hypothesis of \Cref{prop:stratum_H(k^2q)_index_q}. Concretely, we define $\lambda = (k + 1)/2$ if $k \equiv 1 \mod{4}$ and $\lambda = (k + 1)/4$ if $k \equiv 3 \mod{4}$, and assume that every prime factor of $\lambda$ is congruent to $1$ modulo $3$ (in particular, $\lambda$ is odd).
		
		We will exhibit a group $G$ of order $n = 6(k + 1)$, together with two generators $x, y \in G$, such that $[x, y]$ has order $k + 1$. Since the assumed conditions on $\lambda$ are necessary by \Cref{lem:quaternion_klein_x_cyclic}, this is enough to finish the proof.
		
		We know that there exists $r \in \ZZ/\lambda\ZZ$ such that $r^3 = 1$ and such that $r - 1$ has order $\lambda$ by \Cref{prop:derived_sgp_semi_direct_product_vfacile}.
		
		Now, \Cref{lem:quaternion_klein_x_cyclic} suggests considering two cases:
		
		\begin{enumerate}[wide, labelindent=0pt, label=\underline{Case \arabic{*}:}, itemsep=1ex]
			\item Let $G = (\ZZ/\lambda\ZZ \times \ZZ/2\ZZ \times \ZZ/2\ZZ) \rtimes_{\varphi \times \theta} \ZZ/3\ZZ$, where $\varphi \colon \ZZ/3\ZZ \to \Aut(\ZZ/\ZZ\lambda)$ and $\theta \colon \ZZ/3\ZZ \to \Aut(\ZZ/2\ZZ \times \ZZ/2\ZZ)$ are explicitly given by the relations $\varphi(1)(1) = r$, and $\theta(1)(u, v) = (u + v, u)$. Since both $\varphi(1)$ and $\theta(1)$ have order $3$, the semidirect product construction giving $G$ is well-defined.
			
			Consider the elements $x = ((1, 1, 0), 0)$ and $y = ((0, 0, 1), 1)$. We will show that $G = \langle x,y \rangle$ and that $[x, y]$ has order $k + 1 = 2\lambda$.
			
			We compute the commutator:
			\begin{align*}
				[x, y] &= ((1, 1, 0), 0) \cdot ((0, 0, 1), 1) \cdot ((1, 1, 0), 0)^{-1} \cdot  ((0, 0, 1), 1)^{-1} \\
				&= ((1, 1, 0), 0) \cdot ((0, 0, 1), 1) \cdot ((-1, 1, 0), 0) \cdot ((0, 1, 1), -1) \\
				&= ((1, 1, 1), 1) \cdot ((-1, 1, 0), 0) \cdot ((0, 1, 1), -1) \\
				&= ((1 - r, 0, 0), 1) \cdot ((0, 1, 1), -1) \\
				&= ((1 - r, 0, 1), 0).
			\end{align*}
			
			Since $1 - r$ has order $\lambda$ in $\ZZ/\lambda\ZZ$ and $(0, 1)$ has order $2$ in $\ZZ/2\ZZ \times \ZZ/2\ZZ$, which are coprime, we obtain that
			\[
			\ord([x, y]) = \ord(1 - r)\ord(1) = 2\lambda = k + 1.
			\]
			Moreover, since $(1, 0)$ and $(0, 1)$ generate $\ZZ/2\ZZ \times \ZZ/2\ZZ$, and $\lambda$ and $4$ are coprime, \Cref{lem:extension_cyclic} shows that $x' = (1, 1, 0)$ and $y' = (0, 0, 1)$ generate $\ZZ/\lambda\ZZ \times \ZZ/2\ZZ \times \ZZ/2\ZZ$.
			
			Now, observe that $y^2 x^{\lambda} y = (y', 0)$. Indeed,
			\begin{align*}
				y^2 x^{\lambda} y &= ((0, 0, 1), 1) \cdot ((0, 0, 1), 1) \cdot (\lambda(1, 1, 0),0) \cdot ((0, 0, 1), 1) \\
				&= ((0, 1, 1), 2) \cdot ((0, 1, 0), 0) \cdot ((0, 0, 1), 1) \\
				&= ((0, 1, 0), 2) \cdot ((0, 0, 1), 1) \\
				&= ((0, 0, 1), 0) = (y', 0).
			\end{align*}
			
			In particular, $\langle x,y \rangle$ contains $(\ZZ/\ZZ\lambda \times \ZZ/2\ZZ \times \ZZ/2\ZZ) \times \{0\}$. Furthermore, it contains the two elements $y$ and $y^2$, which belong respectively to $(\ZZ/\ZZ\lambda \times \ZZ/2\ZZ \times \ZZ/2\ZZ) \times \{1\}$ and to $(\ZZ/\ZZ\lambda \times \ZZ/2\ZZ \times \ZZ/2\ZZ) \times \{2\}$. Therefore $G = \langle x, y \rangle$.
			\item Let $G = (\ZZ/\lambda\ZZ \times \Q_8) \rtimes_{\varphi \times \theta} \ZZ/3\ZZ$, where $\varphi \colon \ZZ/3\ZZ \to \Aut(\ZZ/\ZZ\lambda)$ and $\theta \colon \ZZ/3\ZZ \to \Aut(\Q_8)$ are explicitly given by the relations $\varphi(1)(1) = r$, and
			\[
			\theta(1)(-1) = -1, \qquad \theta(1)(\qi) = -\qj, \qquad \theta(1)(\qj) = \qk, \qquad \theta(1)(\qk) = -\qi.
			\]
			Since both $\varphi(1)$ and $\theta(1)$ have order $3$, the semidirect product construction giving $G$ is well-defined.
			
			Consider the elements $x = ((1, \qi), 0)$ and $y = ((0, \qk), 1)$. We will show that $G = \langle x, y\rangle$ and that $[x, y]$ has order $k + 1 = 4\lambda$.
			
			We compute the commutator:
			\begin{align*}
				[x, y] &= ((1, \qi), 0) \cdot ((0, \qk), 1) \cdot ((1, \qi), 0)^{-1} \cdot ((0, \qk), 1)^{-1} \\
				&= ((1, \qi), 0) \cdot ((0, \qk), 1) \cdot ((-1, -\qi), 0) \cdot ((0, -\qj), -1) \\
				&= ((1, -\qj), 1) \cdot ((-1, -\qi), 0) \cdot ((0, -\qj), -1) \\
				&= ((1 - r, 1), 1) \cdot ((0, -\qj), -1) \\
				&= ((1 - r, -\qk), 0).
			\end{align*}
			Since $1 - r$ has order $\lambda$ in $\ZZ/\lambda\ZZ$ and $-\qk$ has order $4$ in $\Q_8$, which are coprime, we obtain that
			\[
			\ord([x, y]) = \ord(1 - r)\ord(-\qk) = 4\lambda = k + 1.
			\]
			Moreover, since $\qi$ and $\qk$ generate $\Q_8$, and $\lambda$ and $8$ are coprime, \Cref{lem:extension_cyclic} shows that $x' = (1, \qi)$ and $y' = (0, \qk)$ generate $\ZZ/\lambda\ZZ \times \Q_8$.
			
			Now, choose $\varepsilon \in \{1, 3\}$ so that $\lambda \equiv \varepsilon \mod 4$. We have that $y^2 x^{\varepsilon\lambda} y = (y', 0)$. Indeed,
			\begin{align*}
				y^2x^{\varepsilon\lambda}y & = ((0,\qk),1) \cdot ((0,\qk),1) \cdot ((1,\qi)^{\varepsilon\lambda},0)\cdot ((0,\qk),1)) \\
				& = ((0+\varphi(1)(0),\qk \cdot \theta(1)(\qk)),2) \cdot ((0,\qi),0)\cdot ((0,\qk),1) \\
				&= ((0,-\qj),2) \cdot ((0,\qi),0)\cdot ((0,\qk),1)  \\
				&= ((0,-\qj \cdot \theta(2)(\qi)),2) \cdot ((0,\qk),1)  \\
				&= ((0,\qi),2) \cdot ((0,\qk),1)
				\\
				&=((0,\qi \cdot \theta(2)(\qk)),3) \\
				&=((0,\qk),0) = (y',0).
			\end{align*}
			
			In particular, $\langle x,y \rangle$ contains $(\ZZ/\ZZ\lambda \times \Q_8) \times \{0\}$. Furthermore, it contains the two elements $y$ and $y^2$, which belong respectively to $(\ZZ/\ZZ\lambda \times \Q_8) \times \{1\}$ and to $(\ZZ/\ZZ\lambda \times \Q_8) \times \{2\}$. Therefore $G = \langle x, y \rangle$.
		\end{enumerate}
	\end{proof}
	
	\subsection{Proof of Theorem \ref{thm:no_regular_origamis}} \label{s:no_regular_origamis_proof}
	Assume $g$ is either of the form $g = p+1$ or $g = pq+1$, where $p,q > 3$ are prime numbers. Recall from \Cref{sec:buiding_regular_origamis} that a genus-$g$ regular origami can only belong to a stratum of the form $\calH(m^s)$, where $m \divides (2g - 2)$ and $s = (2g - 2) / m$. Furthermore:
	
	\smallbreak
	\begin{itemize}
		\item since $g-1$ is not divisible by $2$ or $3$, $g \notin \genera(1)$ by the work of Schlage-Puchta and Weitze--Schmithüsen \cite[Theorem 1.1]{Puchta_Schmithusen}, and we can rule out the case $m=1$;
		\item since $g-1$ is not divisible by $2$ or $9$, we can rule out the case $m=2$ by \Cref{thm:H(2^g-1)};
		\item it is well-known that every origami in a minimal stratum $\calH(2g - 2)$ has a trivial automorphism group \cite[Proposition 2.4]{Matheus_Moller_Yoccoz}. In particular, no regular origamis exist in the stratum $\calH(2g - 2)$ and we can rule out the case $m=2g-2$.
	\end{itemize}
	\smallbreak
	
	Now, if $g = p+1$, the only remaining divisor of $2g-2 = 2p$ is $m=p=g-1$, and this case is also ruled out by \Cref{lem:H(g-1 g-1)} since $p$ is odd by assumption. This shows that there are no regular origamis in case \eqref{i:no_regular_origamis_1}.
	
	\medbreak
	
	If $g = pq+1$, we can rule out $m=pq = g-1$ for the same reason. Therefore, we only need to consider the additional divisors $p,q, 2p$, and $2q$. 
	
	The cases $m = p$ and $m = q$ correspond to the symmetric cases $\calH(p^{2q})$ and $\calH(q^{2p})$. Since $p, q > 3$, this case is treated in \Cref{thm:stratum_H(k^2q)}, showing that no regular origamis exist. 
	
	The cases $m = 2p$ and $m = 2q$ correspond to $\calH(2p^q)$ or $\calH(2q^p)$. We know that such strata contain regular origamis by \Cref{cor:Sophie_Germain} if and only if $p=q$ is a Sophie Germain prime. This shows that:
	\begin{itemize}
		\item there are no regular origamis if either \eqref{i:no_regular_origamis_2}, or \eqref{i:no_regular_origamis_3} holds;
		\item if $g = p^2+1$ and $p \geq 5$ is a Sophie Germain prime, every genus-$g$ regular origami belongs to $\calH(2p^p)$ and has translation group $\ZZ/(2p+1)\ZZ \rtimes \ZZ/p\ZZ$. Therefore, $p^2+1 \in \genera (2p)$.
	\end{itemize}
	
	\newpage
	\appendix 
	\section{Explicit computations and summaries}\label{appendix:tables}
	
	The following table shows regular origamis apparently realizing a maximal translation group for each $g$ such that $\translations(g) \leq 2000$. We also include the values of $\orderm(g)$, and of $\constant(g) = \translations(g) / (g - 1) = 2(\orderm(g) + 1) / \orderm(g)$. The examples were found using randomized computer experiments, so some values of $\translations(g)$ could actually be larger, and some regular origamis could be missing.
	
	For all unlisted values of $g$, we have $\translations(g) = 4(g - 1)$ if $g - 1$ is even or divisible by $3$. Otherwise, if no genus-$g$ regular origami exists, then $\translations(g) = 2(g - 1)$ (but, since some regular origamis could be missing from this list, we cannot ensure that this is the case if $g$ is unlisted and $g - 1$ is both odd and not divisible by $3$). We include a translation group with a short description that realizes the maximal number of translations for each genus, although this choice is, in general, not unique.
	
	\begin{table}[H]
		\centering
		\def\arraystretch{1.3}
		\begin{tabular}{|c|c|c|c|c|c|}
			\hline
			$g$ & $\translations(g)$ & $\orderm(g)$ & $\constant(g)$ & Stratum & Translation group \\ \hline
			$26$ & $55$ & $10$ & $11/5$ & $\calH(10^5)$ & $\ZZ/11\ZZ \rtimes \ZZ/5\ZZ$ \\ \hline
			$122$ & $253$ & $22$ & $23/11$ & $\calH(22^{11})$ & $\ZZ/23\ZZ \rtimes \ZZ/11\ZZ$ \\ \hline
			$126$ & $275$ & $10$ & $11/5$ & $\calH(10^{25})$ & $(\ZZ/11\ZZ \rtimes \ZZ/5\ZZ) \times \ZZ/5\ZZ$ \\ \hline
			$176$ & $385$ & $10$ & $11/5$ & $\calH(10^{35})$ & $(\ZZ/11\ZZ \rtimes \ZZ/5\ZZ) \times \ZZ/7\ZZ$  \\ \hline
			$246$ & $497$ & $70$ & $71/35$ & $\calH(70^{7})$ & $\ZZ/71\ZZ \rtimes \ZZ/7\ZZ$  \\ \hline
			$276$ & $660$ & $5$ & $12/5$ & $\calH(5^{110})$ & $\PSL(2, 11)$  \\ \hline
			$326$ & $715$ & $10$ & $11/5$ & $\calH(10^{65})$ & $(\ZZ/11\ZZ \rtimes \ZZ/5\ZZ) \times \ZZ/13\ZZ$ \\ \hline
			$426$ & $935$ & $10$ & $11/5$ & $\calH(10^{85})$ & $(\ZZ/11\ZZ \rtimes \ZZ/5\ZZ) \times \ZZ/17\ZZ$ \\ \hline
			$456$ & $1092$ & $5$ & $12/5$ & $\calH(5^{182})$ & $\PSL(2, 13)$ \\ \hline
			$476$ & $1045$ & $10$ & $11/5$ & $\calH(10^{95})$ & $(\ZZ/11\ZZ \rtimes \ZZ/5\ZZ) \times \ZZ/19\ZZ$ \\ \hline
			$530$ & $1081$ & $46$ & $47/23$ & $\calH(46^{23})$ & $\ZZ/47\ZZ \rtimes \ZZ/23\ZZ$ \\ \hline
			$576$ & $1265$ & $10$ & $11/5$ & $\calH(10^{115})$ & $(\ZZ/11\ZZ \rtimes \ZZ/5\ZZ) \times \ZZ/23\ZZ$ \\ \hline
			$606$ & $1331$ & $10$ & $11/5$ & $\calH(10^{121})$ & $\ZZ/121\ZZ \rtimes \ZZ/11\ZZ$ \\ \hline
			$626$ & $1375$ & $10$ & $11/5$ & $\calH(10^{125})$ & $(\ZZ/11\ZZ \rtimes \ZZ/5\ZZ) \times \ZZ/25\ZZ$ \\ \hline
			$726$ & $1595$ & $10$ & $11/5$ & $\calH(10^{145})$ & $(\ZZ/11\ZZ \rtimes \ZZ/5\ZZ) \times \ZZ/29\ZZ$ \\ \hline
			$776$ & $1705$ & $10$ & $11/5$ & $\calH(10^{155})$ & $(\ZZ/11\ZZ \rtimes \ZZ/5\ZZ) \times \ZZ/31\ZZ$ \\ \hline
			$834$ & $1673$ & $238$ & $239/119$ & $\calH(238^{7})$ & $\ZZ/239\ZZ \rtimes \ZZ/7\ZZ$ \\ \hline
			$842$ & $1711$ & $58$ & $59/29$ & $\calH(58^{29})$ & $\ZZ/59\ZZ \rtimes \ZZ/29\ZZ$ \\ \hline
			$846$ & $1703$ & $130$ & $131/65$ & $\calH(130^{13})$ & $\ZZ/131\ZZ \rtimes \ZZ/13\ZZ$ \\ \hline
			$848$ & $1771$ & $22$ & $23/11$ & $\calH(22^{77})$ & $(\ZZ/23\ZZ \rtimes \ZZ/11\ZZ) \times \ZZ/7\ZZ$ \\ \hline
			$876$ & $1925$ & $10$ & $11/5$ & $\calH(10^{175})$ & $(\ZZ/11\ZZ \rtimes \ZZ/5\ZZ) \times \ZZ/35\ZZ$ \\ \hline
		\end{tabular}
		\label{tab:c(g)}
		\caption{Value of $\translations(g)$ for small $g$.}
	\end{table}
	
	\newpage
	
	For each prime $2 < m \leq 300$ with $m \equiv 2 \mod{3}$, the following table shows the smallest prime $p$ satisfying the conditions in \Cref{lem:infinitely_many_primes}. Namely, we have that $p \equiv \pm 1 \mod{2(m + 1)}$ and that $(p - 1)(p + 1)/(4(m + 1))$ is an integer not divisible by any prime number $q < m$. Hence, the regular origami with translation group $G = \PSL(2, p)$ produces a genus-$g$ surface with $g \in \genera(m)$. We also include $n = |G|$.
	
	\begin{table}[H]
		\centering
		\def\arraystretch{1.2}
		\begin{tabular}{|c|c|c|c|}
			\hline
			$m$ & $p$ & $g$ & $n$ \\ \hline
			$5$ & $11$ & $276$ & $660$ \\ \hline
			$11$ & $23$ & $2784$ & $6072$ \\ \hline
			$17$ & $37$ & $11952$ & $25308$ \\ \hline
			$23$ & $47$ & $24864$ & $51888$ \\ \hline
			$29$ & $59$ & $49620$ & $102660$ \\ \hline
			$41$ & $83$ & $139524$ & $285852$ \\ \hline
			$47$ & $5087$ & $32224176332$ & $65819594208$ \\ \hline
			$53$ & $107$ & $300564$ & $612468$ \\ \hline
			$59$ & $7079$ & $87208034462$ & $177372273480$ \\ \hline
			$71$ & $13967$ & $671699860608$ & $1362320844048$ \\ \hline
			$83$ & $167$ & $1150464$ & $2328648$ \\ \hline
			$89$ & $179$ & $1417860$ & $2867580$ \\ \hline
			$101$ & $23053$ & $3032798528504$ & $6125652473412$ \\ \hline
			$107$ & $24407$ & $3601166766512$ & $7269645061368$ \\ \hline
			$113$ & $227$ & $2898564$ & $5848428$ \\ \hline
			$131$ & $263$ & $4513344$ & $9095592$ \\ \hline
			$137$ & $277$ & $5274912$ & $10626828$ \\ \hline
			$149$ & $44699$ & $22178309490152$ & $44654314409700$ \\ \hline
			$167$ & $56113$ & $43907394117050$ & $88340625289392$ \\ \hline
			$173$ & $347$ & $10385364$ & $20890788$ \\ \hline
			$179$ & $359$ & $11502720$ & $23133960$ \\ \hline
			$191$ & $383$ & $13972224$ & $28090752$ \\ \hline
			$197$ & $397$ & $15563592$ & $31285188$ \\ \hline
			$227$ & $457$ & $23756232$ & $47721768$ \\ \hline
			$233$ & $467$ & $25352964$ & $50923548$ \\ \hline
			$239$ & $479$ & $27360960$ & $54950880$ \\ \hline
			$251$ & $503$ & $31689504$ & $63631512$ \\ \hline
			$257$ & $200723$ & $2013932191752920$ & $4043537007566172$ \\ \hline
			$263$ & $138863$ & $666885785842058$ & $1338842946481392$ \\ \hline
			$269$ & $541$ & $39438360$ & $79169940$ \\ \hline
			$281$ & $563$ & $44455044$ & $89226492$ \\ \hline
			$293$ & $587$ & $50393364$ & $101130708$ \\ \hline
		\end{tabular}
		\caption{Choice of $p$ so $g \in \genera(m)$.}
		\label{tab:c(g)_PSL}
	\end{table}
	
	\newpage
	
	The following table summarizes our current knowledge about the sets $\genera(m)$ for $m \in \{1, \dotsc, 25\}$. As previously stated, $\genera(1)$ was completely classified by Schlage-Puchta and Weitze-Schmithüssen \cite{Puchta_Schmithusen}. Additionally, such sets are empty since $3 \divides m$ or $4 \divides m$ (\Cref{thm:not_regular}). Moreover, if $m \in \{5, 11, 17, 23\}$, then $m$ is prime and satisfies $m \equiv 2 \mod{3}$, so $\genera(m)$ is large (\Cref{thm:G_m_general}). Furthermore, when $m \in \{10, 22\}$, them $m = 2p$ for $p$ a Sophie Germain prime, so $\genera(m)$ is nonempty (\Cref{thm:no_regular_origamis}). Finally, we do not know if $\genera(m)$ is empty or not for the remaining values of $m$.
	
	\begin{table}[H]
		\centering
		\def\arraystretch{1.3}
		\begin{tabular}{|c|c|}
			\hline
			$m$ & $\genera(m)$ \\ \hline
			$1$ & $\{g \geq 2 \st g - 1 \not\equiv \pm 1 \mod{6}\}$ \\ \hline
			$2$ & Empty \\ \hline
			$3$ & Empty \\ \hline
			$4$ & Empty \\ \hline
			$5$ & Contains infinitely long arithmetic progressions \\ \hline
			$6$ & Empty \\ \hline
			$7$ & Empty \\ \hline
			$8$ & Empty \\ \hline
			$9$ & Empty \\ \hline
			$10$ & Nonempty \\ \hline
			$11$ & Contains infinitely long arithmetic progressions  \\ \hline
			$12$ & Empty \\ \hline
			$13$ & ? \\ \hline
			$14$ & ? \\ \hline
			$15$ & Empty \\ \hline
			$16$ & Empty \\ \hline
			$17$ & Contains infinitely long arithmetic progressions  \\ \hline
			$18$ & Empty \\ \hline
			$19$ & ? \\ \hline
			$20$ & Empty \\ \hline
			$21$ & Empty \\ \hline
			$22$ & Nonempty \\ \hline
			$23$ & Contains infinitely long arithmetic progressions \\ \hline
			$24$ & Empty \\ \hline
			$25$ & ? \\ \hline
		\end{tabular}
		\caption{Summary of $\genera(m)$.}
	\end{table}
	
	\newpage
	
	\section{Commutator subgroup of \texorpdfstring{$\ZZ/m\ZZ \rtimes \ZZ/n\ZZ$}{Z/mZ : Z/nZ}}
	In this section, we prove:
	
	\begin{Prop}\label{prop:derived_sgp_semi_direct_product_vfacile}
		Let $u \geq 1$ be an integer and $q$ be a prime number. Then, the following are equivalent:
		\begin{enumerate}
			\item There exists $d \in \{2, \dotsc, u-1\}$ such that the commutator subgroup of the semidirect product $\ZZ/u\ZZ \rtimes_d \ZZ/q\ZZ$ is isomorphic to $\ZZ/u\ZZ$;
			\item There exists $d \in \{2, \dotsc, u-1\}$ with $d^q \equiv 1 \mod{u}$ and $\gcd(d-1,u)=1$;
			\item Every prime factor of $u$ is congruent to $1$ modulo $q$.
		\end{enumerate}
	\end{Prop}
	
	It is clear from the last property that $u$ must be odd for it to hold. Indeed, if $q = 2$, the property states that every prime factor of $u$ is odd. If $q > 2$, the property implies that $2$ is not a prime factor of $u$.
	
	Moreover, the first two properties are equivalent by \Cref{lem:commutator_metacyclic}. We will show that the last two are equivalent. In fact, we will show a slightly more general result which is useful for our purposes: we are interested in the case where a semidirect product realizes the translation group of a regular origami in a given stratum.
	
	We have:
	
	\begin{Prop}\label{prop:derived_sgp_semi_direct_product}
		Let $u, \ell \geq 1$ be integers. Assume $u$ is odd. Then, the following are equivalent:
		\begin{enumerate}[label=(\arabic*)]
			\item \label{i:derived_sgp_1} There exist integers $m,n \geq 1$ and $d \in \{2, \dots, m-1\}$ such that the regular origami induced by the group $\ZZ/m\ZZ \rtimes_d \ZZ/n\ZZ$ and the generators $x=(1,0)$ and $y = (0,1)$ belongs to the stratum $\calH((u-1)^\ell)$;
			\item \label{i:derived_sgp_2} There exist integers $m,n$ and $d \in \{2, \dots, n-1\}$ such that:
			\begin{itemize}
				\item $mn = u \ell$,
				\item $d^n \equiv 1 \mod{m}$,
				\item $\gcd(d-1,m) = m/u$;
			\end{itemize} 
			\item \label{i:derived_sgp_3} For every prime power $p^\alpha$ dividing $u$, we have either
			\begin{itemize}
				\item $p^{\alpha + 1} \divides \ell$, or
				\item $p \equiv 1 \mod q$ for some prime divisor $q$ of $\ell$.
			\end{itemize}
		\end{enumerate}
	\end{Prop}
	
	\begin{Rema}
		If $u$ is even, the number $d$ in the first property exists if and only if $\ell$ is a multiple of $4$. In that case, one can choose $\ZZ/2u\ZZ \rtimes_{-1} \ZZ/(\ell/2)\ZZ$.
	\end{Rema}
	
	\Cref{prop:derived_sgp_semi_direct_product_vfacile} is obtained from \Cref{prop:derived_sgp_semi_direct_product} by taking $m = u$ and $\ell = n = q$.
	We first state the following structure result for cyclic $ p$-groups: 
	
	\begin{Prop}\label{prop:structure_invertibles_p_alpha}
		Let $p \neq 2$ be prime and let $\alpha \geq 1$ be an integer. Then, the group $( \ZZ/p^\alpha \ZZ )^\times$ is cyclic of order $\varphi(p^{\alpha}) = (p - 1)p^{\alpha-1}$. Furthermore, for any $1 \leq \gamma \leq \alpha-1$, the set 
		\[
		H_\gamma = \{ x \in \ZZ/p^{\alpha}\ZZ \st x \equiv 1 \mod {p^{\gamma}} \}
		\]
		is the (unique) subgroup of $( \ZZ/p^{\alpha} \ZZ )^\times$ of order $p^{\alpha-\gamma}$, and is generated by $(1+p^\gamma)$.
	\end{Prop}
	
	\begin{proof}
		The facts that $( \ZZ/p^\alpha \ZZ )^\times$ is cyclic and that $H_1$ is generated by $1 + p$ are well-known \cite[Theorem 6.7]{Rotman_Groups}. The proof readily extends to the general case.
	\end{proof}
	
	We can now prove \Cref{prop:derived_sgp_semi_direct_product}:
	
	\begin{proof}[Proof of \Cref{prop:derived_sgp_semi_direct_product}]
		For any fixed finite group $G$, recall that the stratum $\calH((u - 1)^\ell)$ contains a regular origami with translation group $G$ if and only if there exist generators $x, y \in G$ such that $H = \langle[x,y] \rangle$ has order $u$ and index $\ell$. In the case where $G = \ZZ/m\ZZ \rtimes_d \ZZ/n\ZZ$, this is equivalent to the conditions:
		\begin{enumerate}[label=(\roman*)]
			\item \label{i:proof_B_1} $nm = |G| = u\ell$ (so a subgroup can have order $u$ and index $\ell$);
			\item \label{i:proof_B_2} $d^n \equiv 1 \mod m$ (for the semidirect product to be well-defined); and
			\item \label{i:proof_B_3} $\gcd(d-1, m) = m/u$ (to have $G'\simeq \ZZ/u\ZZ$).
		\end{enumerate}
		Therefore, properties~\ref{i:derived_sgp_1} and \ref{i:derived_sgp_2} are indeed equivalent.
		We will show that properties~\ref{i:derived_sgp_2} and \ref{i:derived_sgp_3} are equivalent.
		
		\paragraph{\textbf{\ref{i:derived_sgp_2} $\implies$ \ref{i:derived_sgp_3}}} Assume the existence of such $m$, $n$ and $d$, and write $m = p_1^{\alpha_1} \dotsm p_r^{\alpha_r}$ and $n = q_1^{\beta_1} \dotsm q_s^{\beta_s}$ for the prime factorization of $m$ and $n$, respectively.
		
		By condition~\ref{i:proof_B_3}, $u$ is a divisor of $m$, so its prime factors are among the $p_1, \dotsc, p_r$. Take $i \in \{1, \dots, r\}$ such that $p_i \divides u$. Observe that $p_i \neq 2$, since $u$ is assumed to be odd. In particular, $( \ZZ/p_i^{\alpha_i} \ZZ )^\times$ is cyclic.
		
		From condition~\ref{i:proof_B_2}, we deduce that $d^n \equiv 1 \mod {p_i^{\alpha_i}}$. In particular, $n$ is a multiple of the order $\ord(d)$ of $d$ in the group $( \ZZ/p_i^{\alpha_i} \ZZ )^\times$. We also deduce that the prime factorization of $\ord(d)$ can be written in terms of the $q_1, \dotsc, q_s$.
		
		Let $\gamma_i \in \NN$ be the multiplicity of $p_i$ in the prime decomposition of $d - 1$. We apply \Cref{prop:structure_invertibles_p_alpha} to the group $(\ZZ/p_i^{\alpha_i}\ZZ)^\times$ and define the groups $H_1, \dotsc, H_{\alpha_i - 1}$ of respective orders $p_i^{\alpha_i - 1}, \dotsc, p_i$.
		
		We consider two disjoint cases:
		
		\begin{enumerate}[wide, labelindent=0pt, label=\underline{Case \arabic{*}:}, itemsep=1ex]
			\item If $\gamma_i = 0$, then $d \not\equiv 1 \mod{p_i}$. This means that $d \notin H_1$. Thus, $d$ is not contained in any of the $H_\delta$, for $\delta \in \{1, \dotsc, \alpha_i - 1\}$, so $\ord(d)$ is not a power of $p_i$. Since $\ord(d)$ divides $|\ZZ/p_i^{\alpha_i}\ZZ| = (p_i - 1)p_i^{\alpha_i - 1}$, we deduce that $\ord(d) = (p_i - 1)p_i^\delta$ for some $\delta \in \{0, \dotsc, \alpha_i - 1\}$. In particular, $(p_i - 1) \divides \ord(d)$. Hence, the prime factorization of $p_i - 1$ can be written in terms of the primes $q_1, \dotsc, q_s$.
			
			Let $j \in \{1, \dotsc, s\}$ such that $q_j \divides (p_i - 1)$. We obtain that $p_i \equiv 1 \mod{q_j}$. By condition~\ref{i:proof_B_1}, $\ell = (m / u) n$ and, by condition~\ref{i:proof_B_3}, $m / u$ is an integer. Thus, $q_j \divides \ell$.
			\item If $\gamma_i > 0$, then $\gamma_i < \alpha_i$. Indeed, $p_i$ divides $u$ by assumption, and we have $u = m/\gcd(m, d- 1)$ by condition~\ref{i:proof_B_3}. We obtain that the multiplicity of $p_i$ in the prime factorization of $u$ is $\alpha_i - \gamma_i$.
			
			We have that $d \equiv 1 \mod{p_i^{\gamma_i}}$, so $d \in H_{\gamma_i}$. We must have $\ord(d) = p_i^{\alpha_i - \gamma_i}$. Indeed, if $\ord(d) = p_i^{\alpha_i - \delta}$ for $\delta > \gamma_i$, then $d \in H_\delta$, so $d \equiv 1 \mod{p_i^\delta}$ and $p_i^\delta \divides (d - 1)$. This contradicts the definition of $\gamma_i$.
			
			Since $\ord(d)$ divides $n$, we deduce that $p_i^{\alpha_i - \gamma_i} \divides n$. Finally, by condition~\ref{i:proof_B_3}, we have that $\ell = m n / u$, so the multiplicity of $p_i$ in the prime factorization of $\ell$ is at least $\alpha_i + (\alpha_i - \gamma_i) - (\alpha_i - \gamma_i) = \alpha_i > \alpha_i - \gamma_i$.
		\end{enumerate}
		
		\paragraph{\textbf{\ref{i:derived_sgp_3} $\implies$ \ref{i:derived_sgp_2}}} We now assume that $u$ and $\ell$ satisfy property~\ref{i:derived_sgp_3}, and we construct $m$, $n$ and $d$ satisfying conditions~\ref{i:proof_B_1}, \ref{i:proof_B_2} and \ref{i:proof_B_3}.
		
		We start by using property~\ref{i:derived_sgp_3} to write the prime factorization of $u$ in a ``split form'' as $u = p_1^{\alpha_1} \dotsm p_r^{\alpha_r} q_1^{\beta_1} \dotsm q_s^{\beta_s}$, where the primes $p_1, \dotsc, p_r$ are congruent to $1$ modulo some prime factor of $\ell$, and $q_j^{\beta_j + 1} \divides \ell$ for each $j \in \{1, \dotsc, s\}$. Since $u$ is odd, none of these prime numbers is $2$.
		
		Moreover, we may write the prime factorization of $\ell$ as $\ell = q_1^{\gamma_1} \cdots q_t^{\gamma_t}$ for some $t \geq s$ and exponents $\gamma_j > \beta_j$ for each $j \in \{1, \dotsc, s\}$. If $i \in \{1, \dotsc, r\}$, we choose $j_i \in \{1, \dotsc, t\}$ such that $p_i \equiv 1 \mod q_{j_i}$.
		
		We define:
		\begin{align*}
			m &= p_1^{\alpha_1} \dotsm p_r^{\alpha_r} q_1^{\gamma_1} \dotsm q_s^{\gamma_s} \\
			n &= q_1^{\beta_1} \dotsm q_s^{\beta_s} q_{s+1}^{\gamma_{s+1}} \dotsm q_t^{\gamma_t}
		\end{align*}
		
		We have $mn = u\ell$, so condition~\ref{i:proof_B_1} is satisfied.
		
		Now, we have $\varphi(p_i^{\alpha_i}) = (p_i - 1) p_i^{\alpha_i - 1}$, so $q_{j_i} \divides \varphi(p_i^{\alpha_i})$ for each $i \in \{1, \dotsc, r\}$. In particular, since the group $(\ZZ/p_i^{\alpha_i}\ZZ)^\times$ is cyclic (as $p_i \neq 2$), \Cref{thm:fundamental_cyclic_groups} ensures the existence of an order-$q_{j_i}$ element $d_i$. By \Cref{prop:structure_invertibles_p_alpha}, we have $d_i \not\equiv 1 \mod{p_i}$ (since, otherwise, $\ord(d_i) \divides p_i^{\alpha_i - 1})$.
		
		By the Chinese remainder theorem, there exists a unique $d \in (\ZZ/m\ZZ)^\times$ with:
		\begin{itemize}
			\item for every $i \in \{1, \dots, r\}$, $d \equiv d_i \mod {p_i^{\alpha_i}}$;
			\item for every $j \in \{1, \dots, s\}$, $d \equiv 1 + q_j^{\gamma_j - \beta_j} \mod {q_j^{\gamma_j}}$.
		\end{itemize}
		
		Observe that $d^{q_{j_i}} \equiv 1 \mod p_i^{\alpha_i}$ for every $i \in \{1, \dotsc, r\}$ by our choice of $d_i$. Moreover, from \Cref{prop:structure_invertibles_p_alpha}, the order of $1 + q_j^{\gamma_j - \beta_j} \in (\ZZ/q_j^{\gamma_j} \ZZ)^\times$ divides $q_j^{\beta_j}$. Therefore, $d^{q_j^{\beta_j}} \equiv 1 \mod q_j^{\gamma_j}$. We get $d^n \equiv 1 \mod m$, yielding condition~\ref{i:proof_B_2}.
		
		Furthermore, using the notation of \Cref{prop:structure_invertibles_p_alpha}, we have $d_i \notin H_\delta$ for every $i \in \{1, \dotsc, r\}$ and $\delta \in \{1, \dotsc, \alpha_i - 1\}$, since, otherwise, its order $q_{i_j}$ would divide a power of $p_i$. This yields:
		\[
		\gcd(d-1,m) = q_1^{\gamma_1 - \beta_1} \dotsm q_s^{\gamma_s - \beta_s} =\frac{p_1^{\alpha_1} \cdots p_r^{\alpha_r} q_1^{\gamma_1} \dotsm q_s^{\gamma_s}}{p_1^{\alpha_1} \dotsm p_r^{\alpha_r} q_1^{\beta_1} \dotsm q_s^{\beta_s}} = \frac{m}{u}.
		\]
		We obtain that condition~\ref{i:proof_B_3} holds, completing the proof.
	\end{proof}
	
	\sloppy\printbibliography
	
\end{document}